\documentclass[12pt]{amsart}
\usepackage{color}
\usepackage{amsmath}
\usepackage{amssymb}
\usepackage{amsthm}
\usepackage{amscd}
\usepackage{eucal} % for \mathcal 
\usepackage{enumitem}
\usepackage{tikz}
\usepackage{hyperref}
\usepackage[utf8]{inputenc}
\usepackage[all]{xy}
\hypersetup{colorlinks=true}
\allowdisplaybreaks

\theoremstyle{plain}
\newtheorem{theorem}{Theorem}[section]

\newtheorem{lemma}[theorem]{Lemma}
\newtheorem{proposition}[theorem]{Proposition}
\newtheorem{corollary}[theorem]{Corollary}

\newtheorem{question}{Question}

\theoremstyle{definition}
\newtheorem{definition}[theorem]{Definition}
\newtheorem{remark}[theorem]{Remark}
\newtheorem{example}[theorem]{Example}

\numberwithin{equation}{section}

\newcommand\fantome[1]{}

\def\Z{\mathbb Z}

\def\cO{\mathcal O}
\def\cL{\mathcal L}

\DeclareMathOperator{\res}{Res}
\DeclareMathOperator{\Log}{Log}
\DeclareMathOperator{\Hom}{Hom}

\DeclareMathOperator{\Frac}{Frac}

\DeclareMathOperator{\GL}{GL}
\DeclareMathOperator{\End}{End}

\DeclareMathOperator{\Lie}{Lie}

\newcommand{\power}[2]{{#1 [[ #2 ]]}}

\newcommand{\F}{\mathbb{F}}
\newcommand{\G}{\mathbb{G}}
\newcommand{\C}{\mathbb{C}}
\newcommand{\cM}{\mathcal{M}}

\newcommand{\bu}{\mathbf{u}}

\newcommand{\bh}{\mathbf{h}}
\newcommand{\bd}{\mathbf{d}}

\newcommand{\bg}{\mathbf{g}}
\newcommand{\be}{\mathbf{e}}

\newcommand{\bz}{\mathbf{z}}

\newcommand{\fs}{\mathfrak{s}}

\newcommand{\inv}{\ensuremath ^{-1}}
\newcommand{\isom}{\ensuremath \cong}

\newcommand{\TT}{\mathbb{T}}
\newcommand{\YY}{\mathbb{Y}}

\newcommand{\N}{\ensuremath \mathbb{N}}

\DeclareMathOperator{\Exp}{Exp}

\DeclareMathOperator{\Mat}{Mat}

\newcommand{\on}{\ensuremath ^{\otimes n}}

\newcommand{\twist}{^{(1)}}

\newcommand{\twistinv}{^{(-1)}}
\newcommand{\twisti}{^{(i)}}
\newcommand{\twistk}[1]{^{(#1)}}

\definecolor{ForestGreen}{rgb}{0.0, 0.5, 0.0}

\author{Nathan Green}

\email{ngreen@latech.edu}

\title[A Motivic Pairing]{A Motivic Pairing and the Mellin Transform in Function Fields}
\begin{document}

\begin{abstract}
We define two pairings relating the $A$-motive with the dual $A$-motive of an abelian Anderson $A$-module. We show that specializations of these pairings give the exponential and logarithm functions of this Anderson $A$-module, and we use these specializations to give precise formulas for the coefficients of the exponential and logarithm functions. We then use these pairings to express the exponential and logarithm functions as evaluations of certain infinite products. As an application of these ideas, we prove an analogue of the Mellin transform formula for the Riemann zeta function in the case of Carlitz zeta values. We also give an example showing how our results apply to Carlitz multiple zeta values.
\end{abstract}

\subjclass[2010]{Primary 11G09; Secondary 11M32, 11M38, 11R58}

\keywords{Mellin transform, Drinfeld modules, Anderson $t$-modules, $t$-motives, multiple zeta values, log-algebraicity}

\date{\today}

\maketitle

\tableofcontents

\section{Introduction}
The Mellin transform is a useful and ubiquitous tool in analytic number theory. It appears in many formulas involving the Riemann zeta function and related special functions. Recall the definition of the Mellin transform for a real-valued function $f(x)$ with suitable decay conditions at $x=0$ and $x=\infty$,
\begin{equation}\label{D:Mellin Def}
\cM(f)(s) = \int_0^\infty f(x)x^{s-1}dx,
\end{equation}
for suitable $s\in \C$ (see \cite{Zag} for details and convergence discussion). A classical formula relates the Mellin transform of the exponential function to the Riemann zeta function and the gamma function. Let $f = 1/(e^x-1)$. Then we have
\begin{equation}\label{E:Mellin Zeta Values}
\cM\left(\frac{1}{e^x-1}\right )(s) = \Gamma(s)\zeta(s).
\end{equation}

One of the main results in this paper establishes an analogue of the above formula for function fields of curves defined over a finite field. Developing integration techniques for characteristic-$p$ function fields has traditionally been very difficult due to the difficulty of defining characteristic-$p$ additive measures. In this paper, we sidestep the need for integration by instead defining maps arising from the motivic structure of certain $t$-modules related to the Carlitz zeta function which take the place of integration in this specific situation.

We now briefly state a special case of the main theorem of this paper and explain why this should be viewed as the correct analogue of the Mellin transform. Let $q=p^r$ be a prime power, and let $A = \F_q[\theta]$ and $K=\F_q(\theta)$ (in most sections of the paper we actually prove our results for an arbitrary function field). Let $K_\infty= \F_q((1/\theta))$ and let $\C_\infty$ be a completion of an algebraic closure of $K_\infty$. Set $u = \frac{\tilde \pi}{\theta-t}\in \C_\infty(t)$, where $\tilde \pi$ is the Carlitz period (see \eqref{E:Carlitz period def}) and $t$ is an independent variable. For a specified element $\bz \in\C_\infty$ (here $n$ is the dimension of a particular $t$-module) we define a map $\delta_{1,\bz}^M$ from a $t$-motive $M$ to $\C_\infty$, which we view as an analogue of integration. We then have the following special case of our main theorem.
\begin{theorem}\label{T:special case of prod} (Special case of Corollary \ref{C:Carlitz zeta prod}) Set $n=1$, $\bz = 1$ and let $M$ be the $t$-motive associated to the Carlitz module $C$. Then we have
\[\delta^M_{1,\bz}\left(\frac{u}{\exp_C(u)}\right ) = \Gamma_A(n)\zeta_A(n) \in K_\infty,\]
where $\exp_C$ is the Carlitz exponential function, and $\Gamma_A$ and $\zeta_A$ are function field versions of the gamma and Riemann zeta function, respectively (see \eqref{E:zeta def} and \eqref{E: Carlitz factorials}).
\end{theorem}

We note that, in the above formula, $\exp_C$ is actually a $t$-linear extension of the Carlitz exponential; we refer the reader to \S \ref{S:Product Formulas} for full details. Similar formulas also hold for integers $n>1$, but the notation is more cumbersome, so we refer the reader to the statement of Corollary \ref{C:Carlitz zeta prod} for details.

We make a brief comparison between the integral formula for the Mellin transform \eqref{D:Mellin Def} and the definition of our map $\delta_{1,\bz}^M$ (see Def. \ref{D:delta^M defs}). The Mellin transform satisfies the functional equation (for suitable $s\in \C$)
\[\cM(xf(x))(s) = \cM(f(x))(s+1),\]
thus it turns multiplication by $x$ against the function into addition in the complex argument. On the other hand, in Proposition \ref{P:A linearity} we show for $\bz \in \C_\infty^n$
\[\delta_{1,\bz}^M(tf(t)) = \delta_{1,\phi_t(\bz)}^M(f(t)),\]
so our map $\delta_{1,\bz}^M$ turns multiplication by $t$ against the function into the \textit{Anderson $A$-module action} in the argument. Since we are in characteristic $p$, it is not natural to expect multiplication to turn into addition (else multiplication by $t^p$ would cause the transform to vanish). The Anderson $A$-module action is the natural replacement.

On the other hand, we observe that the Mellin transform zeta values formula \eqref{E:Mellin Zeta Values} is valid for all complex values $s\neq 1$, whereas our formula Corollary \ref{C:Carlitz zeta prod} is valid only for positive integer values $n$. Thus we interpret our formulas as giving an algebraic interpolation of the Mellin transform formula at integer values, rather than being a true replacement for the Mellin transform. However, strong similarities remain. For example, we are able to recover function field gamma values using $\delta_{1,\bz}^M$ (see Remark \ref{R:Gamma Remark}) which parallels the classical definition of the gamma function as a Mellin transform.

To solidify the analogy between these two Mellin transform formulas, we recall the classical definition of Bernoulli numbers,
\[\frac{x}{e^x-1} = \sum_{n=0}^\infty \frac{B_n x^n}{n!}.\]
Comparing this with the definition of the Carlitz-Bernoulli numbers (originally due to Carlitz \cite{Car35}) as found in \cite{Tha90},
\[\frac{z}{\exp_C(z)} = \sum_{n=0}^\infty \frac{B_{n,C}z^n}{\Gamma(n)},\]
we see that $\frac{z}{\exp_C(z)}$ is the natural function field analogue of the term $\frac{1}{e^x-1}$ (the extra copy of $x$ in the numerator is due to normalization issues). The Bernoulli numbers, both over $\C$ and over function fields, are intimately connect with special values of zeta functions, which strongly supports our assertion that Corollary \ref{C:Carlitz zeta prod} should be viewed as an analogue of the Mellin transform.

The truly exciting part about these new formulas, is that we prove them in a much more general setting than what we have discussed above. We give a ``product formula" (Cor. \ref{C:Exp and Log Products}) for the logarithm and exponential functions of an arbitrary abelian Anderson $A$-module, which in the particular case of tensor powers of the Carlitz module specifies to Theorem \ref{T:special case of prod}. We call it a ``product formula" because it is only a true product formula for simple cases; in general our formulas gives a finite sum of infinite matrix products. In particular, this formula applies to the case of $t$-modules associated to function field multiple zeta values (defined in \eqref{D:MZVdef}). As an example application of our main theorems, we give a formula of the form
\[\text{``product formula"} = (\theta^2+\theta) \zeta_A(1,3),\]
where the left hand side is a sum of terms which look like a higher dimensional version of $\frac{u}{\exp_C(u)}$, expressed as an infinite product (see \eqref{E:MZV product formula}). Thus, we are able to realize function field multiple zeta values (MZVs) using a generalization of our Mellin transform formula. This mirrors the situation for complex valued MZVs, which can be expressed as a higher-dimensional Mellin transform (see \cite[(2.1)]{FKMT17}). In Section \ref{S:Product Formulas} we give one specific example of this behavior; we plan to study such formulas in more generality in a future project. In Remark \ref{R:MZV remark} we discuss how zeta values potentially appear in the terms of our ``product formula" and we are hopeful that the LHS of our formula can be decomposed into terms involving zeta values and MZV, thus giving new relations between MZV. This also is the topic of a future project.

The major innovation in this paper is a new technique we develop, which we call a motivic pairing. This pairing allows us to give explicit formulas for the coefficients of the exponential and logarithm functions of Anderson $A$-modules (see Def. \ref{D:Anderson A-module}; original definition in \cite{And86}). Our formulas apply in a very general setting, and give large improvements over previously known formulas for the coefficients of the exponential and logarithm. For example, they recover the powerful results contained in \cite[Prop. 2.2]{ANDTR20a} in the case of $A$-finite $t$-modules.

We now briefly discuss the setting of our main theorems and briefly describe the contents of each section. An Anderson $A$-module $\phi$ of dimension $d$ is an $\F_q$-algebra homomorphism 
\[\phi:A \to \Mat_d(\C_\infty)[\tau],\]
where $\Mat_d(\C_\infty)[\tau]$ is a skew polynomial ring where we view $\tau$ acting via the $q$-power Frobenius (see \ref{D:skew ring}). Associated to the Anderson $A$-module $\phi$, there is an exponential function $\Exp_\phi:\C_\infty^d \to \C_\infty^d$ which serves a similar role as the exponential map of a Lie group. The inverse function of $\Exp_\phi$ is called the logarithm function, $\Log_\phi$. Viewed as functions from $\C_\infty^d\to \C_\infty^d$, both $\Exp_\phi$ and $\Log_\phi$ are $\F_q$-linear power series with coefficients in $\Mat_d(\C_\infty)$. For each Anderson $A$-module we functorially define an associated $A$-motive and dual $A$-motive. We review this theory in detail in Section \ref{S:Review}.

The main results of Sections \ref{S:Exp Pairing} and \ref{S:Log Pairing} give explicit formulas for the coefficients of the exponential and logarithm functions of an arbitrary Anderson $A$-module subject to a couple technical conditions (see Defs. \ref{D:A motive} and \ref{D:Xi Regular}). These formulas are contained in Corollaries \ref{C:Exp coeff} and \ref{C:Log Coeff}. Our main technique is to construct two pairings between the associated $A$-motive and dual $A$-motive of $\phi$ (see Def. \ref{D:A motive} for the definition of $A$-motives). We then show that a specialization of these parings returns the exponential and logarithmic functions and this allows us to prove precise formulas for $\Exp_\phi$ and $\Log_\phi$ in Theorems \ref{T:Exponential from Pairing} and \ref{T:Log from Pairing}. We view this $A$-motivic set up as being a very natural setting to study these coefficients, and our formulas subsume many previous ad hoc approaches (see \cite{Pap}, \cite{Gre17b}, \cite{EGP13}). As an immediate application of our formulas, we give a sufficient condition for the exponential and logarithm coefficients of $t$-modules to be invertible in Theorems \ref{T:Invertibility of P_i} and \ref{T:Invertibility of Q_i}.

Using this motivic set up, in Section \ref{S:Product Formulas} we develop a pair of ``product formulas" for the exponential and logarithm functions in Theorem \ref{T:Product Formulas}. We specialize these product formulas to the case of $t$-modules related to zeta values and multiple zeta values to give a function field version of the Mellin transform formula for the Riemann zeta function in Theorem \ref{C:Carlitz zeta prod} and Example \ref{Ex:MZV (3,1)}. We also give an additional brief application of these new techniques in Section \ref{S:Log Alg} to recover and slightly extend a result of Anderson related to log-algebraicity.

\section{Review of Anderson Motives and Modules}\label{S:Review}
In this section we review the definitions and much of the theory involved in $A$-motives, dual $A$-motives and Anderson $A$-modules. Most of these results are originally due to Anderson \cite{And86}, although here we follow the exposition of Hartl and Juschka \cite{HJ20}. We also prove several new results which will be used in subsequent sections to define our motivic pairings. Let $\F_q$ be the finite field of size $q=p^r$ and let $X$ be a projective, smooth, geometically connected curve defined over $\F_q$ with fixed point at infinity $\infty \in X(\F_q)$. Let $K$ be the function field of $X$ and let $A$ be the functions on $X$ regular away from $\infty$. Let $K_\infty$ be the completion of $K$ at the place of $K$ corresponding to $\infty$, and let $\C_\infty$ be the completion of an algebraic closure of $K_\infty$. Let $L$ be an algebraically closed subfield with $K\subset L\subset \C_\infty$. We comment that in some subsections we will set $L=\C_\infty$, and the reader doesn't lose much generality by simply making this assumption throughout the paper. We view $L\otimes_{\F_q} A$ as the ring of functions regular away from the point $\infty$ of the curve $L\times_{\F_q} X$ and view $\Frac(L\otimes_{\F_q} A)$ as its function field. Finally, we let $\Xi$ be the point of $L\times_{\F_q} X$ corresponding to the kernel of the multiplication map, $L\otimes_{\F_q} A \to L$ given by $\sum l_i \otimes a_i \mapsto \sum l_i a_i$.

\begin{definition}\label{D:skew ring}
Let $(L\otimes_{\F_q} A)[\tau]$ be the skew polynomial ring with commutativity relation
\[\tau (\ell \otimes a) = (\ell^q \otimes a) \tau,\quad \forall \ell\otimes a \in  L\otimes_{\F_q} A,\]
where everything else commutes. We similarly define $(L\otimes_{\F_q} A)[\sigma]$, where $\sigma = \tau\inv$. When there is no confusion, we will sometimes refer to $L\otimes 1 \subset L \otimes A$ simply as $L \subset L\otimes A$ and similarly for $A \subset L\otimes A$.
\end{definition}

\begin{definition}
Let $J$ be the maximal ideal of $L\otimes_{\F_q} A$ generated by $(1\otimes a)-(a\otimes 1)$ for all $a\in A$ (recall $A\subset L$).
\end{definition}

\begin{definition}\label{D:A motive}
An abelian $A$-motive of rank $r$ and dimension $d$ is an $(L\otimes_{\F_q} A)[\tau]$-module $M$ which is a free $L[\tau]$ module (via $L\otimes 1$) of rank $d$ and a finitely generated, projective $L\otimes A$ module of constant local rank $r$ such that for $n\in \Z$ sufficiently large we have
\[J^nM\subset \tau M.\]
Note: Here $\tau M$ is a $L\otimes A$-submodule of $M$ because we have assumed $L$ is algebraically closed.
\end{definition}

\begin{definition}\label{D:Dual A Motive}
An $A$-finite dual $A$-motive of rank $r$ and dimension $d$ is an $(L\otimes_{\F_q} A)[\sigma]$-module $N$ which is a free $L[\sigma]$ module (via $L\otimes 1$) of rank $d$ and a finitely generated, projective $L\otimes A$ module of constant local rank $r$ such that for $n\in \Z$ sufficiently large we have
\[J^nN\subset \sigma N.\]
\end{definition}

\begin{remark}
In this article, we will assume all the $A$-motives and dual $A$-motives are Abelian and $A$-finite, respectively, and will hence shorten this to simply $A$-motive or dual $A$-motive.
\end{remark}

\begin{definition}\label{D:Frobenius def}
For a simple tensor $h = z\otimes f \in L\otimes A$, we define the $i$th Frobenius twist of $h$ to be
\[h\twisti = z^{q^i}\otimes f \in L\otimes A,\]
then extend linearly to sums of simple tensors. We extend this definition coordinate-wise to vectors and matrices of elements in $L\otimes A$. We also extend this definition to points on $X$ by twisting each coordinate of $X$.
\end{definition}

\begin{definition}\label{D:Anderson A-module}
A $d$-dimensional Anderson $A$-module ($E$,$\phi$) is an $A$-module scheme $E$ over $L$ satsifying the following two properties:
\begin{enumerate}
\item As an $\F_q$-vector space scheme over $L$, we have $E$ is isomorphic to $\G_a^d$.
\item Under the identification $\End_{\F_q}(E) = \End_{\F_q}(\G_a^d) = \Mat_d(L)[\tau]$ the action of $a \in A$, seen as an $\F_q$-linear vector space scheme endomorhpism $\phi_a$, has the form
\[\phi_a = d[a] + A_1 \tau + \dots ,\quad A_i \in \Mat_d(L)\]
where $d[a] = a I + N$ for some nilpotent matrix $N\in \Mat_d(L)$ (depending on $a$). 
\end{enumerate}
A Drinfeld $A$-module is a dimension $d=1$ Anderson $A$-module which has positive degree in $\tau$. A $t$-module is a $d$-dimensional Anderson $A$-module for which $A=\F_q[t]$ (we have $X=\mathbb P^1$).

The map $d:A\to \Mat_n(\overline K_\infty)$, where $d[a]$ is the constant term of $\phi_a$, is a ring homomorphism, and we will use this notation regularly throughout the paper. We will refer to $\F_q$-vector space scheme over $L$ which this $A$-action generates as $\Lie_E$. For more details on these constructions we refer the reader to \cite[\S 5]{Gos96}
\end{definition}

\begin{definition}\label{D:Staroperator}
We define a ring homomorphism $*:(L\otimes A)[\tau] \to (L\otimes A)[\sigma]$ by setting, for $\alpha = \sum_{i=0}^m a_i \tau^i\in (L\otimes A)[\tau]$,
\[\alpha^* = \sum_{i=0}^m a_i\twistk{-i} \sigma^i.\]
We extend this map to matrices $B\in \Mat_{m\times n}((L\otimes A)[\tau])$ by setting $B^*$ to be matrix $B^\top$ with the $*$ map applied to each coordinate (here $^\top$ is the transpose). Note that for two such matrices we have $(B_1B_2)^* = B_2^*B_1^*$.
\end{definition}

To each Anderson $A$-module $\phi$, there exists an associated Anderson $A$-motive and dual $A$-motive, which we denote by $M_\phi$ and $N_\phi$. This correspondence is given functorially by
\[M_\phi = \Hom_{\F_q}(E,\G_a)\isom \Mat_{1\times d}(L[\tau]),\quad N_\phi = \Hom_{\F_q}(\G_a,E)\isom L[\tau]^d,\]
(the underlying algebraic group of $\phi$ is $\G_a^d$) where $a\in A$, $b\in L$ and $\tau$ act on $m\in M_\phi$ by
\[a \cdot m = m\circ \phi_a,\quad b\cdot m = bm,\quad \tau\cdot m = \tau m,\]
(here $\tau$ acts as the $q$-power Frobenius) and $a\in A$, $b\in L$ and $\sigma$ act on $n\in N_\phi$ by
\[a \cdot n = \phi_a\circ n,\quad b\cdot n = bn,\quad \sigma\cdot m = m \tau.\]
Thus $M_\phi$ and $N_\phi$ are $(L\otimes A)[\tau]-$ and $(L\otimes A)[\sigma]-$modules, and it is shown in \cite[\S 2.5.2]{HJ20} that they are $A$-motives and dual $A$-motives, respectively.  Note that we may identify $N_\phi = \Mat_{1\times d}(L[\sigma])$, in which case the actions can be written as
\[a \cdot n =  n \circ\phi_a^*,\quad b\cdot n = bn,\quad \sigma\cdot m = \sigma m.\]
For details and proofs of the above facts we refer the reader to \cite[\S 2.5]{HJ20} and \cite[\S 4-5]{BP20}. In a recent paper \cite{Mau21}, Maurischat proved an important relationship between $M_\phi$ and $N_\phi$.

\begin{theorem}[Maurischat]
For an Anderson $A$-module $\phi$, we have that $M_\phi$ is abelian if and only if $N_\phi$ is $A$-finite.
\end{theorem}

Each Anderson $A$-module $\phi$ has an associated exponential and logarithm function. The exponential function is the unique $\F_q$-linear power series, denoted
\begin{equation}\label{D:Exp and Log}
\Exp_\phi(\bz) = \sum_{i=0}^\infty Q_i \bz\twisti,
\end{equation}
for $\bz \in \C_\infty^d$ with $Q_i \in \Mat_{d\times d}(\C_\infty)$, such that $Q_0 = I$ and
\[\Exp_\phi(d[a]\bz) = \phi_a(\Exp_\phi(\bz))\]
for all $a\in A$. The exponential function converges for any $\bz\in \C_\infty^d$

The logarithm function is defined as the formal power series inverse of $\Exp_\phi$ and is denoted
\[\Log_\phi(\bz) = \sum_{i=0}^\infty P_i \bz\twisti.\]
It is the unique $\F_q$-linear function with $P_i=I$ satisfying
\[\Log_\phi(\phi_a\bz) = d[a](\Log_\phi(\bz).\]
We note that the logarithm function, in contrast with the exponential, has some bounded domain of convergence in $\C_\infty^d$.

Throughout the remainder of this paper we will fix an Anderson $A$-module $\phi$. We will denote the associated $A$-motive and dual $A$-motive by $M$ and $N$ respectively, suppressing the $\phi$ from their notation.

Let $M_K = M\otimes_{L\otimes_{\F_q} A} \Frac(L\otimes_{\F_q} A)$ be the $(\Frac(L\otimes_{\F_q} A))[\tau]$-module where the $\tau$-action on $M$ is as given above and the $\tau$-action on $\Frac(L\otimes_{\F_q} A)$ acts as the Frobenius (as in Def. \ref{D:Frobenius def}). Define $N_K = N\otimes_{L\otimes_{\F_q} A} \Frac(L\otimes_{\F_q} A)$ similarly where the $\sigma$ acts on $\Frac(L\otimes_{\F_q} A)$ as the inverse of Frobenius. Note that $M_K$ and $N_K$ are $\Frac(L\otimes_{\F_q} A)$-vector spaces of dimension $r$.

\begin{lemma}\label{L:Bases}
Let $\phi$ be an Anderson $A$-module and let $M$ and $N$ be the associated $A$-motive and dual $A$-motive, respectively.
\begin{enumerate}
\item There exist $\Frac(L\otimes_{\F_q} A)$-bases for $M_K$ and $N_K$, denoted $\{c_1,\dots,c_r\}$ and $\{d_1,\dots,d_r\}$ respectively, such that if the $\tau$-action on $M_K$ in this basis is given by
\[\tau (c_1,\dots,c_r)^\top = \Phi (c_1,\dots,c_r)^\top,\]
for $\Phi \in \Mat_{r\times r}(\Frac(L\otimes_{\F_q} A))$, then the $\sigma$-action on $N$ is given by
\[\sigma (d_1,\dots,d_r)^\top = \Phi^\top (d_1,\dots,d_r)^\top.\]
\item There exists an $L[\tau]$-basis for $M$ and an $L[\sigma]$-basis for $N$, denoted $\{g_1,\dots,g_d\}$ and $\{h_1,\dots,h_d\}$ respectively, such that for $b \in 1\otimes A$, if the action of $b$ on $M$ is given by the matrix
\[b (g_1,\dots,g_d)^\top = \Theta_b (g_1,\dots,g_d)^\top,\]
for $\Theta_b \in \Mat_{d\times d}(L[\tau])$, then the action of $b$ on $N$ is give by the matrix
\[b (h_1,\dots,h_d)^\top = \Theta_b^* (h_1,\dots,h_d)^\top,\]
where $\Theta_b^*$ is defined in \ref{D:Staroperator}. We will define $\bg = (g_1,\dots,g_d)$ and $\bh = (h_1,\dots, h_d)$.
\end{enumerate}
\end{lemma}
\begin{proof}
The proof of (1) follows from \cite[Prop. 2.4.3]{HJ20}. The proof of part (2) follows from \cite[(5.2)]{Mau21} (see also \cite[\S 4.4]{BP20}).
\end{proof}

\begin{definition}\label{D:Xi Regular}
We fix the bases of Lemma \ref{L:Bases} for use throughout the paper. Further, given such bases and the matrix $\Phi$ as defined above, we say that the Anderson $A$-module $\phi$ is $\Xi$-regular if each coordinate of $\Phi$ is regular at $\Xi\twisti$ for all $i\in \Z$. Note that if $A=\F_q[\theta]$, then $\Phi \in \Mat_d(L[t])$ and $\det(\Phi) \in (t-\theta)^dL[t]$, thus all $t$-modules are automatically $\Xi$-regular. Throughout the rest of the paper we will assume that all our Anderson $A$-modules are $\Xi$-regular.
\end{definition}

\begin{remark}\label{R:Remark on Xi regularity}
We also comment here, thanks to excellent suggestions by the referee, that many of the constructions in this paper can be done without assuming $\
\Xi$-regularity. Notably the construction of the pairings $F$ and $G$ of Sections \ref{S:Exp Pairing} and \ref{S:Log Pairing} can be done without it, as is explained in Remark \ref{R:Basis free remark}. In fact, it is a widely assumed folklore conjecture that all Anderson $A$-modules are $\Xi$-regular, so we hope that with future work this assumption can be removed entirely.
\end{remark}

\begin{definition}\label{D:delta^N defs}
For an Anderson $A$-module $\phi$, let $N$ be the associated dual $A$-motive with rank $r$ and dimension $d$ and let $\{h_1,\dots,h_d\} \subset N$ be an $L[\sigma]$-basis for $N$ as described in Lemma \ref{L:Bases}. We define two maps,
\[\delta_0^N,\, \delta_1^N:N\to L^d,\]
for $n\in N$ by first expressing $n$ in the basis $\{h_1,\dots,h_d\}$,
\begin{align}\label{D:N coeff basis}
\begin{split}
n  =& c_{1,0}h_1 + c_{1,1}\sigma(h_1) + \dots + c_{1,m}\sigma^m(h_1) \\
&+ c_{2,0}h_2 + c_{2,1}\sigma(h_2) + \dots + c_{2,m}\sigma^m(h_2)\\
&\quad\vdots\\
&+ c_{d,0}h_d + c_{d,1}\sigma(h_d) + \dots + c_{d,m}\sigma^m(h_d),
\end{split}
\end{align}
(some of the $c_{j,m}$ could be zero) then writing
\begin{equation}
\delta_0^N(n) = 
\left (\begin{matrix}
c_{1,0}\\
c_{2,0}\\
\vdots\\
c_{d,0}
\end{matrix}\right )
\end{equation}
and
\begin{equation}\label{D:delta_n^1 Def}
\delta_1^N(n) = 
\left (\begin{matrix}
c_{1,0}\\
c_{2,0}\\
\vdots\\
c_{d,0}
\end{matrix}\right ) +
\left (\begin{matrix}
c_{1,1}\\
c_{2,1}\\
\vdots\\
c_{d,1}
\end{matrix}\right )\twist+\dots + 
\left (\begin{matrix}
c_{1,m}\\
c_{2,m}\\
\vdots\\
c_{d,m}
\end{matrix}\right )\twistk{m}.
\end{equation}
\end{definition}

\begin{definition}\label{D:delta^M defs}
For an Anderson $A$-module $\phi$, let $M$ be the associated $A$-motive with rank $r$ and dimension $d$ and let $\{g_1,\dots,g_d\} \subset M$ be an $L[\tau]$-basis for $M$ as described in Lemma \ref{L:Bases}. Fix an element $\bz \in L^d$. We define two maps,
\[\delta_0^M:M\to L^d,\quad \delta_{1,\bz}^M:M\to L\]
for $m\in M$ by first expressing $m$ in the basis $\{g_1,\dots,g_d\}$,
\begin{align}\label{D:M coeff basis}
\begin{split}
m  =& c_{1,0}g_1 + c_{1,1}\tau(g_1) + \dots + c_{1,m}\tau^m(g_1) \\
&+ c_{2,0}g_2 + c_{2,1}\tau(g_2) + \dots + c_{2,m}\tau^m(g_2)\\
&\quad\vdots\\
&+ c_{d,0}g_n + c_{d,1}\tau(g_d) + \dots + c_{d,m}\tau^m(g_d),
\end{split}
\end{align}
then writing
\begin{equation}
\delta_0^M(m) = 
\left (\begin{matrix}
c_{1,0}\\
c_{2,0}\\
\vdots\\
c_{d,0}\\
\end{matrix}\right ).
\end{equation}
For the definition of $\delta_{1,\bz}^M$, we recall that $M = \Hom_{\F_q}(\phi,\mathbb G_a) \isom \Mat_{1\times d}(L[\tau])$, where $a\in A$ and $b\in L$ act on $m\in M$ via
\[a m = m\circ \phi_a,\quad bm = b\cdot m.\]
We then define
\begin{equation}
\delta_{1,\bz}^M(m) = m(\bz).
\end{equation}
Note that this definition depends on the element $\bz \in L^d$, but in our applications we will often use it for an arbitrary element of $L^d$, so when convenient we will suppress this dependency in our notation and simply write $\delta_1^M$.

We comment that $\delta_0^M$ and $\delta_0^N$ can be recovered by the maps induced by $M \to M/\tau M$ and $N\to N/\sigma N$. Similarly, $\delta_1^N$ can be recovered by the map induced by $N \to N/(1-\sigma)N$, and while it is true that $\delta_{1,\bz}^M$ can be recovered by the map induced by $N\to N/(1-\tau)N$, this is much more subtle.
\end{definition}

\begin{proposition}\label{P:A linearity}
For $\phi$, $m\in M$, and $n\in N$ as given above, we have the following for all $a\in A$.
\begin{enumerate}
\item $\delta_1^N(an) = \phi_a(\delta_1^N(n))$
\item $\delta_{1,\bz}^M(am) = \delta_{1,\phi_a(\bz)}^M(m)$
\item $\delta_0^N(an) = d[a]\delta_0^N(n)$
\item $\delta_0^M(am) = d[a]^\top \delta_0^M(m)$. 
\end{enumerate}
\end{proposition}
\begin{proof}
These statements all follow from \cite[Prop. 2.4.3 and Prop. 2.5.8]{HJ20}. They can also easily be seen by direct calculation following from the definitions.
\end{proof}

\begin{definition}
For $c\in \C_\infty$, we define the Tate algebra in the variable $t$ with radius $|c|_\infty$, where $|\cdot|_\infty$ is the norm used for the completion $K_\infty$. 
\[\TT_c = \biggl\{ \sum_{i=0}^\infty b_i t^i \in \power{\C_\infty}{t} \biggm| \big\lvert c^i b_i \big\rvert \to 0 \biggr\}.\]
We then identify $\F_q[t]$ as a subring of $A$ via a fixed injective ring morphism $\F_q[t]\to A$ (this is not unique)  and define the ring
\begin{equation}\label{D:Tate algebras}
\YY_c = (\C_\infty\otimes_{\F_q} A)\otimes_{\C_\infty[t]} \TT_c.
\end{equation}
Note that this ring is a ring of regular functions on an open rigid analytic subspace of the rigid analytic space associated to our curve $X$. We refer the reader to \cite[\S 2.3.3]{HJ20} for full details on the rigid analytic geometry (in the case $d_\infty=1$).
% In the same vein, we also define $\Gamma'$ to be the ring of rigid analytic functions on the affine curve $(\C_\infty \times_{\F_q} X)\setminus \{\infty\}$ which are regular away from $\Xi\twisti$ for all $i\geq 1$ (we note for the reader's sake that $\Gamma' = \mathcal O(\dot{\mathfrak C}_\C\setminus \cup_{i\geq 1} V(\sigma^{i*}J))$ in the notation of \cite{HJ20}). 
We let $\tau$ (and $\sigma$) act on $\YY_c$ diagonally via Frobenius twisting (inverse twisting). Thus we let $L=\C_\infty$ and set $M_\theta = M\otimes_{(\C_\infty\otimes A)} \YY_\theta$ and $N_\theta = N\otimes_{(\C_\infty\otimes A)} \YY_\theta$ with diagonal $\tau$ and $\sigma$ action, respectively. We set $N_\Gamma$ to be the submodule of $N_\theta$ consisting of all $h\in N_\theta$ such that $\sigma^i(h)\in N_\theta$ for all $i \in \Z$.

We note here that whenever we are dealing with a Tate algebra, we will set $L=\C_\infty$, so that our coefficients live in a complete field. However, this will mostly occur in Section \ref{S:Log Alg}, so most of the paper will be developed for general $L$.
\end{definition}

\begin{definition}\label{D:TT norms}
We define a norm on $\TT_c$ for $h = \sum_{i=0}^\infty b_i t^i\in \TT_c$, by setting
\[\lVert h \rVert_c = \sup |c|_\infty^i \cdot |b_i|_\infty.\]
We then extend this norm to $\YY_c$ by recognizing that $\YY_c \isom \TT_c^\ell$ where $\ell$ is the degree of $K/\F_q(t)$, and taking the max of the norm in each coordinate. Finally, we extend this norm to matrices $\Mat_{m\times m}(\YY_c)$, again by taking the max of the norm of each coordinate. Note that this matrix norm is sub-multiplicative, so for $B_1,B_2 \in \Mat_{m\times m}(\YY_c)$ we have
\[\lVert B_1 B_2\rVert_c \leq \lVert B_1 \rVert_c \cdot \lVert B_2\rVert_c.\]
Finally, we comment that the ring $\YY_c$ is complete in the defined norm. In particular in this paper, we are usually concerned with the norm $\lVert \cdot \rVert_\theta$, where $\theta \in \C_\infty$ an image of $t$ under the inclusion $\F_q[t]\to A\to \C_\infty$.
\end{definition}

\begin{proposition}\label{P:delta_0 extensions}
There exists an extension of the map $\delta_0^N$ to a $\C_\infty\otimes A$-module homomorphism
\[\delta_0^N:N_\theta \to \C_\infty^d.\]
Let $S\subset L\otimes A$ be the set of elements not in $J$ and let $(L\otimes A)_{(J)}$ be the localization of $L\otimes A$ at $S$. Then there exists an extension of $\delta_0^M$ to $M_{(J)} = M\otimes (L\otimes A)_{(J)}$ (here $M$ is defined over an arbitrary $L$).
\end{proposition}

\begin{proof}
The proof of the first extension of $\delta_0^N$ is given in \cite[Prop 2.5.8]{HJ20} and relies on the fact that the map $\delta_0^N$ factors through $N/J^dN$ ($J$ is the maximal ideal generated by $(1\otimes a)-(a\otimes 1)$). For the second, we note that $M/\tau M$ is a $J$-torsion $L\otimes A$-module, so $\delta_0^M M\to M/\tau M \isom \Hom_L(\Lie_\phi(L),L)$ factors uniquely as
\[
\xymatrix{
M \ar[d]\ar[r]^{\delta_0^M}  & M/\tau M \\
M_{(J)}\ar@{-->}[ur] &
}
\]
where the up-right direction arrow represents the extension of $\delta_0^M$ described in the proposition.
\end{proof}

\begin{definition}\label{D:N[[sigma]]}
Note that $N\otimes_{L[\sigma]} L[[\sigma]]\isom L[[\sigma]]^d$, and thus for $n\in  N\otimes_{L[\sigma]} L[[\sigma]]$, similarly to \eqref{D:N coeff basis}, we may write for $c_{i,j}\in L$
\begin{align}\label{D:N coeff infinite}
\begin{split}
n  =& c_{1,0}h_1 + c_{1,1}\sigma(h_1) + \dots  \\
&+ c_{2,0}h_2 + c_{2,1}\sigma(h_2) + \dots  \\
&\quad\vdots\\
&+ c_{d,0}h_n + c_{d,1}\sigma(h_d) + \dots  .
\end{split}
\end{align}
We let $N\langle \sigma \rangle $ denote the subset of $n\in N\otimes_{L[\sigma]} L[[\sigma]]$ such that we have $(c_{1,j},\dots,c_{d,j})\twistk{j}\to 0$ as $j\to \infty$ in the above expansion for $n$. Thus there exists an extension of the map $\delta_1^N:N\langle \sigma \rangle \to \C_\infty^d$ defined by applying $\delta_1^N$ to partial sums of \eqref{D:N coeff infinite}, then taking the limit. For the $A$-motive $M$ and for fixed $\bz\in L^d$, we define $M\langle \tau \rangle $ similarly, except subject to the condition that $(c_{1,j},\dots,c_{d,j})\bz\twistk{j}\to 0$. As above, there exists an extension of $\delta_{1,\bz}^M$ to $M\langle \tau \rangle $.
\end{definition}

Using the $\Frac(L\otimes_{\F_q} A)$-basis $\{c_1,\dots, c_r\}$ from Lemma \ref{L:Bases}, if we view $M_K \isom (\Frac(L\otimes_{\F_q} A))^r$, then for $(a_1,\dots,a_r) \in \Mat_{1\times r}(\Frac(L\otimes_{\F_q} A))$, we find that
\[\tau (a_1,\dots,a_r) (c_1,\dots,c_r)^\top =  (a_1,\dots,a_r)\twist \tau (c_1,\dots,c_r)^\top =  (a_1,\dots,a_r)\twist \Phi (c_1,\dots,c_r)^\top.\]
After transposing, we may write the $\tau$ action of $M_K$ viewed as $(\Frac(L\otimes_{\F_q} A))^r$ as
\[\tau (a_1,\dots,a_r)^\top = \Phi^\top((a_1,\dots,a_r)\twist)^\top .\]
Similarly we find that $N_K \isom (\Frac(L\otimes_{\F_q} A))^r$ and
\[\sigma (a_1,\dots,a_r)^\top = \Phi((a_1,\dots,a_r)\twistinv)^\top .\]

\begin{definition}\label{D:tau inverse}
For $m \in  M_K$ we use the isomorphism in the preceding discussion to write $m=(a_1,\dots,a_r)\in  (\Frac(L\otimes_{\F_q} A))^r$ and define
\[\tau\inv(m) = ((\Phi\inv)^\top)\twistinv((a_1,\dots,a_r)\twistinv)^\top.\]
Note that the coordinates of $\Phi\inv$ are regular at $\Xi\twisti$ for $i\neq 0$ by Definition \ref{D:Xi Regular} and by the Cayley-Hamilton inverse formula, thus the coordinates of $(\Phi\inv)\twisti$ are regular at $\Xi$ for all $i<0$. Thus if $m \in M \subset M_K$ then $\tau\inv(m)\in M_{(J)}$ and, then by Proposition \ref{P:delta_0 extensions} we may evaluate $\delta_0^M$ at $\tau^{-i}(m)$. We similarly define for $n\in N_K$
\[\sigma\inv(n) = (\Phi\inv)\twist((a_1,\dots,a_r)\twist)^\top.\]
Note that if $n \in N\subset N_K$, then any poles of $\sigma\inv(n)$ will be at $\Xi\twisti$ for $i\geq 1$, and thus we may evaluate $\delta_0^N$ at $\sigma\inv(n)$. We further extend this $\sigma$-action to $n \in N_\theta$ and by similar arguments conclude that $\sigma\inv(n)$ is in $N_\theta$ and thus may be evaluated under $\delta_0^N$.
\end{definition}

We note that $\tau(\tau\inv(m)) = \tau\inv(\tau(m)) = m$ and $\sigma(\sigma\inv(n)) = \sigma\inv(\sigma(n)) = n$ for all $m\in M_K$ and $n\in N_K$.

\section{The Exponential Motivic Pairing}\label{S:Exp Pairing}
The goal of this section is to define a pairing between the rings $(L\otimes_{\F_q} A)[\tau]$ and $(L\otimes_{\F_q} A)[\sigma]$, as well as a related pairing between the motives $M$ and $N$, which specializes to give the exponential function. This will allow us to give very explicit formulas for the coefficients of the exponential function for any abelian $\Xi$-finite Anderson $A$-module.

\begin{definition}
Let $\phi$ be a $\Xi$-regular Anderson $A$-module. Let $M$ and $N$ be the associated $A$-motive and dual $A$-motive. Let $\{g_1,\dots,g_d\}$ and $\{h_1,\dots,h_d\}$ be the bases described in Lemma \ref{L:Bases} above. For $x\in (L\otimes_{\F_q} A)[\tau]$ and $y \in (L\otimes_{\F_q} A)[\sigma]$ and $\bz \in L^d$, define the pairing
\[F:(L\otimes_{\F_q} A)[\tau]\times (L\otimes_{\F_q} A)[\sigma] \to \C_\infty^d,\]
\begin{equation}\label{D:Fpairing}
F(x,y;\bz) = \sum_{i=0}^\infty \sum_{k = 1}^{d} \delta_1^N\left(\delta_0^M(\tau^{-i}(x(g_k)))^\top \bz \sigma^i(y(h_k)))  \right )
\end{equation}
\end{definition}

\begin{remark}
We recall for the reader's convenience that $\delta_0^M(g)\in L^d$ for any $g\in M_{(J)}$. Thus, in the above summand, $\delta_0^M(\tau^{-i}(x(g_k)))^\top \bz \in L$, and thus it makes sense to multiply this by $\sigma^i(y(h_k))\in N$. This is an element of $N$, so we can evaluate it under $\delta_1^N$, and taking the sum for $i\geq 0$ gives an element $\Mat_d(L)[[\bz]]$. The question remains whether this sum converges to an element in $\C_\infty^d$. For now, we view the pairing as a formal sum in $\Mat_d(L)[[\bz]]$, and we will address the convergence issues in Lemma \ref{L:F convergence} after establishing some of the properties of $F$.
\end{remark}

It will be useful to have a finer version of the pairing $F$. To this end, for any two elements $g\in M$ and $h\in N$ we define
\begin{equation}\label{D:Fpairinganything}
F(x,y;\bz;g,h) = \sum_{i=0}^\infty \delta_1^N\left(\delta_0^M(\tau^{-i}(x(g)))^\top \bz \sigma^i(y(h)))  \right ).
\end{equation}
The pairing $F(x,y;\bz;g,h)$ exists and converges defined based on the same reasoning as Lemma \ref{L:F convergence}.

\begin{remark}\label{R:Basis free remark}
We thank the referee for the following helpful suggestion. It is possible to give a completely coordinate free definition of \eqref{D:Fpairinganything}. Namely, for $m\in M$ and $n\in N$ and $\bz\in \Lie_E(L)$ we set
\[F(m,n;\bz) = \sum_{i=0}^\infty \delta_1^N\left (\delta_0^M(\tau^{-i}m)(z)\cdot \sigma^i n \right )\in E(L),\]
where we comment that this definition makes sense because we identify $M/\tau M \isom \Hom_L(\Lie_E(L),L)$ as in \cite[2.5.2]{HJ20}, thus $\delta_0^M(\tau^{-i}m)(z)$ makes sense and is in $L$. Then, $N$ is an $L$-vector space, so $\delta_0^M(\tau^{-i}m)(z)\cdot \sigma^i n \in N$. Finally, $\delta_1^N\left (\delta_0^M(\tau^{-i}m)(z)\cdot \sigma^i n \right )\in N/(1-\sigma)N \isom E(L)$. Then, to recover Definition \ref{D:Fpairing} we evaluate the above function $F$ at a chosen basis for $M$ and its dual basis for $N$ and take the trace.

We comment that, $\tau^{-i}m\in M_{(J)}$ always, with no restrictions on $M$. So with this construction, the $\Xi$-regularity condition of Definition \ref{D:Xi Regular} is not needed. Only once we choose bases for $M$ and $N$ as $\Frac(L\otimes_{\F_q} A)$-vector spaces, is the $\Xi$-regularity condition needed to evaluate this representation of $\tau^{-i}m$ at $\Xi$.

\end{remark}

\begin{remark}
We comment that $F(x,y;\bz)$ should be viewed as a trace of $F(x,y;\bz,g,h)$ over the basis of Lemma \ref{L:Bases}. Namely,
\[F(x,y;\bz) = \text{Trace}\left( \left [\sum_{i=0}^\infty \delta_1^N\left(\delta_0^M(\tau^{-i}(x(g_j)))^\top \bz \sigma^i(y(h_k)))  \right )\right ]_{j,k = 1}^{d}\right ). \]
\end{remark}

\begin{lemma}\label{L:Fine bilinearity}
For any $g\in M$ and $h\in N$, the pairing $F(x,y;\bz;g,h)$ satisfies:
\begin{enumerate}
\item $F(ax,y;\bz;g,h) = F(x,ay;\bz;g,h)$ for all $a\in L$.
\item $F(\tau x,y;\bz;g,h) = F(x,\sigma y;\bz;g,h)$.
\item $F(bx,y;\bz;g,h) = F(x,y;d[b]\bz;g,h)$ for all $b\in A$.
\item $F(x,by;\bz;g,h) = \phi_b\circ F(x,y;\bz;g,h)$ for all $b\in A$.
\end{enumerate}
\end{lemma}

\begin{proof}
The first equality of part (1) follows by observing that $\tau\inv a = a^{1/q} \tau\inv$, that $\sigma a = a^{1/q} \sigma$ and recalling that $\delta_0^M$ is $L$-linear. Statement (2) follows from the observation that $\delta_0^M(\tau(g)) = 0$ for any $g \in M$. Parts (3) and (4) follow from the $A$-linearity of $\delta_0^M$ and $\delta_1^N$ as described in Proposition \ref{P:A linearity}.
\end{proof}

\begin{proposition}\label{P:F Bilinearity}
For the pairing $F(x,y;\bz)$, we have
\begin{enumerate}
\item $F(a x,y;\bz) = F(x,ay;\bz)$ for all $a\in L$.
\item $F(a x,y;\bz) = F(x,ay;\bz) = aF(x,y;\bz)$ for all $a\in L$, if $x,y\in L$.
\item $F(\tau x,y;\bz) = F(x,\sigma y;\bz)$.
\item $F(bx,y;\bz) = F(x,by;\bz)= \phi_b\circ F(x,y;\bz)= F(x,y;d[b]\bz)$ for all $b \in A$ and $x,y\in A$.
\end{enumerate}
\end{proposition}

\begin{remark}
Before giving the proof of Proposition \ref{P:F Bilinearity}, we make a brief comment about calling $F$ a pairing, which will apply to the other pairings presented in the paper. The above proposition does not actually show that $F$ is an $(L\otimes A)[\tau]$-bilinear pairing. It would rather be more accurate to describe $F$ as a mapping which is simultaneously $A$-bilinear, $L$-bilinear and which is $\tau$- and $\sigma$-symmetric. That being said, we will call $F$ a pairing for ease of discussion.
\end{remark}

\begin{proof}
The proof of parts (1) and (3) follow directly from Lemma \ref{L:Fine bilinearity}. The proof of part (2) follows from (1) together with the fact that $\delta_1^N(a\sigma^i(h_j)) = a^{q^i}\delta_1^N(\sigma^i(h_j))$ for the basis elements $h_j$. For the proof of part (4), we fundamentally use the fact that the action of $A$ on $M$ and on $N$ is related via the $^*$ operator of Definition \ref{D:Staroperator}. Recall from Lemma \ref{L:Bases} that there exists a matrix $\Theta_b \in \Mat_{d\times d}(L[\tau])$ such that for if we label $\bg = (g_1,\dots,g_d)^\top$ and $\bh = (h_1,\dots,h_d)^\top$
\[b \bg = \Theta_b \bg,\text{ and}\quad b \bh = \Theta_b^* \bh\]
where the elements of $\bg$ form an $L[\tau]$-basis for $M$ and $\bh$ an $L[\sigma]$-basis for $N$. Let us label the $j$th row of $\Theta_b$ as $\Theta_{b,j}$. Then, using the fact that $x\in  A$ and thus commutes with $b$ and $g_k$ we get,
\begin{align*}
F(bx,y;\bz) &= \sum_{i=0}^\infty \sum_{k = 1}^{d} \delta_1^N\left(\delta_0^M(\tau^{-i}(bx(g_k)))^\top \bz \sigma^i(y(h_k))  \right )\\
&= \sum_{i=0}^\infty \sum_{k = 1}^{d} \delta_1^N\left(\delta_0^M(\tau^{-i}( x(bg_k)))^\top \bz \sigma^i(y(h_k))  \right )\\
&= \sum_{i=0}^\infty \sum_{k = 1}^{d} \delta_1^N\left(\delta_0^M(\tau^{-i}( x(\Theta_{b,k}\bg)))^\top \bz \sigma^i(y(h_k))  \right ).
\end{align*}
Then, using Lemma \ref{L:Fine bilinearity} parts (1) and (2) on the individual terms of $\Theta_{b,k}\bg$, we get
\begin{align*}
\sum_{i=0}^\infty \sum_{k = 1}^{d} \delta_1^N\left(\delta_0^M(\tau^{-i}( x(\Theta_{b,k}\bg)))^\top \bz \sigma^i(y(h_k))  \right ) &= \sum_{i=0}^\infty \sum_{k = 1}^{d} \delta_1^N\left(\delta_0^M(\tau^{-i}( x(g_k)))^\top \bz \sigma^i(y(\Theta_{b,k}^*\bh))  \right )\\
&= \sum_{i=0}^\infty \sum_{k = 1}^{d} \delta_1^N\left(\delta_0^M(\tau^{-i}( x(g_k)))^\top \bz \sigma^i(by(h_k))  \right )\\
& = F(x,by;\bz).
\end{align*}
The final two equalities of part (4) follow from Lemma \ref{L:Fine bilinearity} parts (3) and (4).
\end{proof}

\begin{theorem}\label{T:Exponential from Pairing}
For $\bz \in L^d$, as formal series in $\Mat_d(L)[[\bz]]$ we have
\[F(1,1;\bz) = \Exp_\phi(\bz).\]
\end{theorem}

\begin{proof}
Observe first that $F(1,1;\bz)$ may be expressed as an $\F_q$-linear power series in $\bz$. Further, its first nonzero term is $\text{Id}_d \bz$. It remains to show that $F(1,1;\bz)$ satisfies the correct functional equation given in Def. \ref{D:Exp and Log}. Using Proposition \ref{P:F Bilinearity}(3), we see that for all $a\in A$
\begin{align*}
F(1,1;d[a]\bz) &= \delta_1^N\left(\sum_{i=0}^\infty \sum_{k = 1}^{d} \delta_0^M(\tau^{-i}(g_k))^\top d[a]\bz \sigma^i(h_k)  \right )\\
&= \delta_1^N\left(\sum_{i=0}^\infty \sum_{k = 1}^{d} (d[a]^\top\delta_0^M(\tau^{-i}(g_k)))^\top \bz \sigma^i(h_k)  \right )\\
&= \delta_1^N\left(\sum_{i=0}^\infty \sum_{k = 1}^{d} \delta_0^M(\tau^{-i}(ag_k))^\top \bz \sigma^i(h_k)  \right )\\
&= \delta_1^N\left(\sum_{i=0}^\infty \sum_{k = 1}^{d} \delta_0^M(\tau^{-i}(g_k))^\top \bz \sigma^i(ah_k)  \right )\\
&= \delta_1^N\left(a\sum_{i=0}^\infty \sum_{k = 1}^{d} \delta_0^M(\tau^{-i}(g_k))^\top \bz \sigma^i(h_k)  \right )\\
&= \phi_a \left(\delta_1^N\left(\sum_{i=0}^\infty \sum_{k = 1}^{d} \delta_0^M(\tau^{-i}(g_k))^\top \bz \sigma^i(h_k)  \right) \right )\\
&=\phi_a(F(1,1;\bz)).
\end{align*}
Thus $F(1,1;\bz)$ equals the exponential function for $\phi$.
\end{proof}

\begin{corollary}\label{C:Exp coeff}
Let $\phi$ be a $\Xi$-regular Anderson $A$-module (\ref{D:Xi Regular}). If we write
\[\Exp_\phi(\bz) = \sum_{i=0}^\infty Q_i \bz\twisti,\]
then
\[Q_i = (\delta_0^M(\tau^{-i}(g_1))\twisti,\dots,\delta_0^M(\tau^{-i}(g_d))\twisti)^\top.\]
\end{corollary}

\begin{remark}\label{R:Choice of basis}
We comment briefly on the dependence of our formulas in Theorem \ref{T:Exponential from Pairing} and Corollary \ref{C:Exp coeff} on the basis chosen in Lemma \ref{L:Bases}. Once we have chosen $L[\tau]$- and $L[\sigma]$-bases $\bg$ and $\bh$ satisfying the conditions in Lemma \ref{L:Bases}, given any matrix $P\in \GL_d(L[\tau])$ we may switch to the basis $\bg_1 = P\bg$. This has the effect for any $b\in A$ that
\[b\bg_1 = (P\Theta_b P\inv) \bg_1.\]
This amounts to changing the coordinates of $\phi$ to an isomorphic Anderson $A$-module given by $P\phi P\inv$ (see \cite[\S 3.5]{NP21}). We then make a corresponding change of basis for $\bh$ given by $\bh_1 = (P^*)\inv\bh$, from which we deduce that
\[b\bh_1 = ((P^*)\inv \Theta_b^* P^*)\bh_1 = (P\Theta_b P\inv)^* \bh_1.\]
Thus, these new bases also satisfy the conditions of Lemma \ref{L:Bases}. In this sense, we are free to choose any $L[\tau]$-basis we want for $M$, then the choice of $L[\sigma]$-basis is forced for $N$. Indeed, under the identification $M=\Hom(E,\G_a)$ and $N = \Hom(\G_a,E)$ if we choose a basis $\bg$ for $M$, then the required basis for $N$ is simply the basis dual to $\bg$ in the sense of $\F_q$-vector space group schemes. So, changing $L[\tau]$-basis for $M$ produces an isomorphic Anderson $A$-module which has exponential function $P\cdot \Exp_\phi(P\inv\bz)$ - as can be seen from our formulas after substituting $\bg_1$ for $\bg$.

On the other hand, changing $\Frac(L\otimes_{\F_q} A)$-bases for $M_K$ and $N_K$ has no effect on the Anderson $A$-module $\phi$ (again, see \cite[\S 3.5]{NP21})), and thus has no effect on our formulas for the exponential. However, it does affect how we calculate $\tau^{-i}(g_k)$ as well as the evaluation of the $\delta_0^M$ map (albeit in opposite ways so that our overall formulas are not affected).
\end{remark}

\begin{proof}
This follows immediately from Theorem \ref{T:Exponential from Pairing} after noting that $\delta_1^N \left( \sigma^i (h_k) \right) = \be_{k}$, the $k$th standard basis vector.
\end{proof}

\begin{lemma}\label{L:F convergence}
The infinite sum
\[\sum_{i=0}^\infty \sum_{k = 1}^{d} \delta_0^M(\tau^{-i}(x(g_k)))^\top \bz \sigma^i(y(h_k)))  \in N\langle \sigma \rangle ,\]
($N\langle \sigma \rangle $ is defined in \ref{D:N[[sigma]]}) thus we may actually write
\[F(x,y;\bz) = \delta_1^N\left(\sum_{i=0}^\infty \sum_{k = 1}^{d} \delta_0^M(\tau^{-i}(x(g_k)))^\top \bz \sigma^i(y(h_k)))  \right ).\]
\end{lemma}

\begin{proof}
Since $F(1,1;\bz)$ and $\Exp_\phi(\bz)$ agree as formal power series, and since we know $\Exp_\phi(\bz)$ is an entire function on $\C_\infty^d$ (see \cite{And86}), we conclude that $F(1,1;\bz)$ is also an entire function of $\bz$. By Proposition \ref{P:F Bilinearity} we see that $F(cx,dx;\bz) = cdF(x,y;\bz)$ for $c,d\in L$, that $F(ax,bx;\bz) =\phi_b \circ F(x,y;d[a]\bz)$ for $a,b\in A$. Further, recall that for $c\in \C_\infty$ and $h_k$ a basis element of Lemma \ref{L:Bases}, if $\delta_1^N(c \sigma^i(h_k))=w\in \C_\infty^d$, then by \eqref{D:delta_n^1 Def} we have that $\delta_1^N(c\sigma^{i+1}(h_k)) = w\twist$. From this, it follows that $F(\tau ,1;\bz) = F(1,\sigma ;\bz) = F(1,1;\bz)\twist$. Thus,  the convergence of $F(x,y;\bz)$ reduces to the convergence of $F(1,1;\bz)$, which is an entire function. Thus $\sum_{i=0}^\infty \sum_{k = 1}^{d} \delta_0^M(\tau^{-i}(x(g_k)))^\top \bz \sigma^i(y(h_k)))  \in N\langle \sigma \rangle $.

\end{proof}

\begin{remark}
We note that the pairing $F:(L\otimes A)[\tau]\times (L\otimes A)[\sigma] \to \C_\infty^d$ is not perfect. Indeed, $M$ and $N$ are torsion $(L\otimes A)[\tau]-$ and $(L\otimes A)[\sigma]-$ modules, respectively, and one can use this torsion to show that the induced map
\[(L\otimes A)[\sigma] \to \Hom_{(L\otimes A)[\tau]}((L\otimes A)[\tau],\C_\infty^d)\]
is not injective. However, if we restrict the pairing to $F:L\times L \to \C_\infty^d$, then the induced map is injective as long as $\Exp_\phi(\bz)\neq 0$. Further, if one restricts the pairing to $A\times A \to \C_\infty^d$, then a short calculation shows that the induced map is actually the Anderson $A$-module morphism $\phi$.
\end{remark}

\begin{example}[Tensor Powers of Carlitz]\label{E:Tensor Power of Carlitz}
In this example we give details of many of the above constructions in the case of the $n$th tensor power of the Carlitz module for fixed $n>0$. We have $K=\F_q(t)$ and $A = \F_q[t]$ and we set $L=\C_\infty$. Let us identify $1\otimes A = \F_q[t]$ and $A\otimes 1 = \F_q[\theta]$ for independent variables $t$ and $\theta$, thus obviating the need for the tensor product in this situation (thus $\C_\infty$ is written using the $\theta$ variable). We may realize the $t$-motive $M$ and the dual $t$-motive $N$ as $M=N=\C_\infty[t]$ (here $L\otimes_{\F_q} A = \C_\infty\otimes_{\F_q} A = \C_\infty[t]$) with
\[\tau(m) = (t-\theta)^n m\twist,\quad \sigma(s) = (t-\theta)^n s\twistinv,\]
for $m\in M$ and $s\in N$. The associated Anderson $A$-module (in this case called a $t$-module) is denoted $C^{\otimes n}$ and is given by
\begin{equation}\label{E:Carlitz Def}
C^{\otimes n}_t = \begin{pmatrix}
\theta & 1 & 0 & \cdots & 0\\
0 & \theta & 1 & \cdots & 0\\
\vdots &\vdots &\vdots  & & \vdots \\
0 & 0 & 0 & \cdots & \theta
\end{pmatrix}
+
\begin{pmatrix}
0 & \cdots & 0\\
\vdots & & \vdots \\
\tau & \cdots & 0
\end{pmatrix}
\end{equation}
(we refer the reader to \cite[\S 3.3]{BP20} for details). In this case, $M$ and $N$ both have the same $\C_\infty[\tau]-$ and $\C_\infty[\sigma]-$ bases, namely
\[\{1,(t-\theta),\dots,(t-\theta)^{n-1}\},\]
and in order to satisfy the conditions of Lemma \ref{L:Bases} we set $g_k = (t-\theta)^{k-1}$ and $h_k= (t-\theta)^{n-k}$ for $1\leq k\leq n$. Further, we have $M_K = N_K = \C_\infty(t)$ and
\[\tau\inv(m) = \frac{1}{(t-\theta^{1/q})^n} m\twistinv\in \C_\infty(t)\]
\[\sigma\inv(m) = \frac{1}{(t-\theta^{q})^n} m\twist\in \C_\infty(t).\]
Let us denote $D_i(t) = (t-\theta)\dots(t-\theta^{q^{i-1}})$, so that
\[(\tau^{-i}(1)\big|_{t=\theta})^{q^i} = \frac{1}{D_i(\theta^{q^i})^n}.\]
The maps $\delta_0^M$ and $\delta_0^N$ may each be given by hyperderivatives. Namely, for $m\in M$, we write
\[m = m(\theta) + \partial_t(m)\big|_{t=\theta} (t-\theta) + \partial^2_t(m)\big|_{t=\theta} (t-\theta)^2 + \dots,\]
where $\partial^k_t(m)$ denotes the $k$th partial hyperderivative of $m$ with respect to the variable $t$. Because $\delta_0^M$ factors through $M/(t-\theta)^nM$, we find
\[\delta_0^M(m) = (m(\theta),\partial_t(m))\big|_{t=\theta},\dots,\partial^{n-1}_t(m)\big|_{t=\theta})^\top.\]
Similarly, we calculate for $s\in N$
\[\delta_0^N(s) = (\partial^{n-1}_t(s)\big|_{t=\theta},\dots,\partial_t(s))\big|_{t=\theta},s(\theta))^\top.\]
Thus, using the product rule for hyperderivatives, we find that
\begin{equation}\label{E:product rule}
\delta_0^M(m)^\top \delta_0^N(n) = \partial_t^n(mn)\big|_{t=\theta}.
\end{equation}
We refer the reader to \cite[\S 3.2-3.4]{CGM19} and \cite[\S 2.3-2.4]{Pap} for full details on hyperderivatives. Our Corollary \ref{C:Exp coeff} then shows that the coefficients of the exponential are given by
\[Q_i = 
\begin{pmatrix}
\frac{1}{D_i(t)^n} & \partial_t(\frac{1}{D_i(t)^n}) & \dots & \partial_t^{n-1}(\frac{1}{D_i(t)^n})\\
\frac{t-\theta}{D_i(t)^n} & \partial_t(\frac{t-\theta}{D_i(t)^n}) & \dots & \partial_t^{n-1}(\frac{t-\theta}{D_i(t)^n})\\
\vdots & \vdots & \ddots & \vdots\\
\frac{(t-\theta)^{n-1}}{D_i(t)^n} & \partial_t(\frac{(t-\theta)^{n-1}}{D_i(t)^n}) & \dots & \partial_t^{n-1}(\frac{(t-\theta)^{n-1}}{D_i(t)^n})
\end{pmatrix}\Bigg|_{t=\theta^{q^i}}
.\]
This recovers Proposition 4.3.6(b) of \cite{Pap}. 
\end{example}

\begin{example}
In this example we give some details showing how our formulas recover results of Thakur \cite[Prop 0.3.6]{Tha93}. We let $X$ be a projective, smooth, geometically connected curve defined over $\F_q$. Thakur shows that we may define $M = \Gamma(U, \cO_{E}(V))  = \bigcup_{i \geq 0} \cL((V) + i(\infty)),$ where $V$ is a particular divisor called the Drinfeld divisor. We endow $M$ with a $\tau$-action given by $\tau z = f z\twist$, where $f$ is the shtuka function (see \cite[(16)]{GP18} or \cite[0.3]{Tha93} for details on these constructions). One can show that $M$ is a free rank 1 $L[\tau]$-module with basis $\{1\}$ and a projective (locally free) $L\otimes A$-module. A short computation shows that we have
\[\tau\inv z = \frac{1}{f\twistinv}z\twistinv,\]
for $z\in M_K$. Finally, Corollary \ref{C:Exp coeff} (here $d=1$) shows that the exponential coefficients $Q_i$ are given by
\[Q_i = \delta_0^M(\tau^{-i}(1)) =  \left (\frac{1}{f\twistinv \cdots f\twistk{-i}}\bigg|_{\Xi}\right )\twisti = \frac{1}{f \cdots f\twistk{i-1}}\bigg|_{\Xi\twisti},\]
which matches Thakur's \cite[Prop 0.3.6]{Tha93}. Similar constructions can be made to recover the author's \cite[Thm. 4.1]{Gre17b}.
\end{example}

\section{The Logarithmic Motivic Pairing}\label{S:Log Pairing}
In this section we will define a pairing similar to $F$ in the previous section which will give the logarithm function upon specialization. This will again give us an extremely explicit formula for the coefficients of the logarithm function. Many of the proofs in this section follow nearly identically to those in Section \ref{S:Exp Pairing}, so we leave many of the details to the reader. Our Theorem \ref{T:Log from Pairing} and Corollary \ref{C:Log Coeff} generalizes formulas for the logarithm function of certain $t$-modules given by Angl\`es, Ngo Dac and Ribeiro in \cite[Prop. 2.2]{ANDTR20a}. In this section, we also give a sufficient criteria for the coefficients of the exponential and logarithm function to be invertible.

\begin{definition}\label{D:G Pairing}
Recall that $\Log_\phi$ has some bounded domain of convergence $D\subset\C_\infty^d$. Let $\bz \in D\cap L^d$ be fixed. Also recall that $\{g_i\}$ and $\{h_i\}$ are the bases of Lemma \ref{L:Bases}. We define a second pairing 
\[G:(L\otimes A)[\sigma]\times (L\otimes A)[\tau] \to \C_\infty^d,\]
for $x\in (L\otimes A)[\sigma]$ and $y\in (L\otimes A)[\tau]$ by setting
\[G(x,y;\bz) = \sum_{i=0}^\infty \sum_{k = 1}^{d} \delta_0^N\left(\sigma^{-i}(x(h_k)) \delta_{1,\bz}^M(\tau^i(y(g_k)))  \right ).\]
For now we view $G(x,y;\bz)$ as a formal $\F_q$-linear power series in the variable $\bz$ (this follows quickly from the definition of $\delta_{1,\bz}^M$ \eqref{D:delta^M defs}). However, we will see in Theorem \ref{T:Log from Pairing} that for the given $\bz$ the pairing coincides with the logarithm function, and thus it actually converges.
\end{definition}

As in the previous section, we also define a finer pairing which will be useful in our proofs. For any $g\in M$ and $h\in N$ we set
\begin{equation}\label{D:G pairing finer}
G(x,y;\bz;g,h) = \sum_{i=0}^\infty \sum_{k = 1}^{d} \delta_0^N\left(\sigma^{-i}(x(h)) \delta_{1,\bz}^M(\tau^i(y(g)))  \right ).
\end{equation}

\begin{lemma}\label{L:G Fine Bilinearity}
For any $g\in M$ and $h\in N$, the pairing $G(x,y;\bz;g,h)$ satisfies:
\begin{enumerate}
\item For all $a\in L$, we have $G(ax,y;\bz;g,h) = G(x,ay;\bz;g,h)$.
\item We have $G(\sigma x,y;\bz;g,h) = G(x,\tau y;\bz;g,h)$.
\end{enumerate}
\end{lemma}

\begin{proof}
The proof is similar to Lemma \ref{L:Fine bilinearity}, and we leave the details to the reader.
\end{proof}

\begin{proposition}\label{P:G Bilinearity}
For the pairing $G(x,y;\bz)$, we have
\begin{enumerate}
\item $G(a x,y;\bz) = G(x,ay;\bz)$ for all $a\in L$.
\item $G(\sigma x,y;\bz) = G(\tau x,y;\bz)$
\item $G(bx,y;\bz) = G(x,by;\bz)$ for all $b \in A$ and $x,y\in A$.
\end{enumerate}
\end{proposition}

\begin{proof}
This proposition follows nearly identically to Proposition \ref{P:F Bilinearity}. Again, we leave the details to the reader.
\end{proof}

\begin{theorem}\label{T:Log from Pairing}
Recall that $\Log_\phi$ has some bounded domain of convergence $D\subset\C_\infty^d$. Let $\bz \in D\cap L^d$ be fixed. We have
\[G(1,1;\bz) = \Log_\phi(\bz).\]
\end{theorem}

\begin{proof}
As in the proof of Theorem \ref{T:Exponential from Pairing}, we observe first that $G(1,1;\bz)$ is a formal $\F_q$-linear power series in $\bz$ with linear term $\bz$. It remains to show that $d[a]G(1,1;\bz) = G(1,1;\phi_z(\bz))$. We do this using Proposition \ref{P:G Bilinearity}(3) and give an abbreviated version of the calculation here,
\begin{align*}
G(1,1;\phi_a(\bz)) &= \delta_0^N\left(\sum_{i=0}^\infty \sum_{k = 1}^{d} \sigma^{-i}(h_k) \delta_{1,\bz}^M(\tau^i(g_k(\phi_a)))  \right )\\
&= \delta_0^N\left(\sum_{i=0}^\infty \sum_{k = 1}^{d} \sigma^{-i}(h_k) \delta_{1,\bz}^M(\tau^i(t g_k))  \right )\\
&= \delta_0^N\left(\sum_{i=0}^\infty \sum_{k = 1}^{d} \sigma^{-i}(th_k) \delta_{1,\bz}^M(\tau^i(g_k))  \right )\\
&= d[a]\delta_0^N\left(\sum_{i=0}^\infty \sum_{k = 1}^{d} \sigma^{-i}(h_k) \delta_{1,\bz}^M(\tau^i(g_k))  \right )\\
&= d[a] G(1,1;\bz).
\end{align*}
\end{proof}

\begin{corollary}\label{C:Log Coeff}
Let $\phi$ be a $\Xi$-regular Anderson $A$-module (\ref{D:Xi Regular}). If we write
\[\Log_\phi(\bz) = \sum_{i=0}^\infty P_i \bz\twisti,\]
then
\[P_i = \left(\delta_0^N(\sigma^{-i}(h_1)),\dots,\delta_0^M(\sigma^{-i}(h_d))\right ).\]
\end{corollary}

\begin{proof}
This follows immediately from Theorem \ref{T:Log from Pairing} by noting that if we write $\bz = (z_1,\dots,z_d)^\top$, then $\delta_{1,\bz}^M(\tau^i(g_k)) = z_k^{q^i}$.
\end{proof}

\begin{remark}
We comment that in the case of the $n$th tensor power of the Carlitz module, the above formula recovers Proposition 4.3.6(a) of \cite{Pap}. 
\end{remark}

\subsection{Invertibility of Exponential and Logarithm Coefficients}

Our formulas in Corollaries \ref{C:Exp coeff} and \ref{C:Log Coeff} also allow us to give a sufficient condition for the coefficients of the logarithm and exponential functions to be invertible in the case of $t$-modules (when $A = \F_q[t]$) for $L=\C_\infty$. This problem was first considered for Drinfeld $A$-modules by Thakur in \cite[3.2]{Tha93}, where he gave two sufficient conditions on the curve $X$ for when the exponential coefficients are invertible (nonzero). Since our formulas depend somewhat on the choice of basis (see Remark \ref{R:Choice of basis}), our goal in this subsection are to give a sufficient condition under which it is possible to find coordinates for the $t$-module in which the coefficients of the exponential are invertible. Our starting point is \cite[Prop. 3.5.7]{NP21}, which gives an explicit formula for the map $\delta_0^N$ in terms of hyperderivatives. This proposition assumes that the dual $t$-motive $N$ (and corresponding $t$-module) has a certain structure. However, in \cite[3.5.11]{NP21} the authors prove that after a suitable change of basis, any dual $t$-motive can be made to satisfy these condition, thus it is really not a restriction. The conditions on the dual $t$-motive $N$ are given as follows.
\begin{enumerate}
\item There exists $C\in \text{GL}_r(\C_\infty[t])$ such that the $\sigma$-action on $N$ is represented by a matrix
\[\Phi = C \begin{pmatrix}
(t-\theta)^{\ell_1} & & \\
 & \ddots &\\
  & & (t-\theta)^{\ell_r}
\end{pmatrix},\]
such that $\ell_k \geq 0$ with $\sum \ell_k = d$, the $\C_\infty[\sigma]$-dimension of $N$.
\item If $e_i$ are basis elements for $N$ as a $\C_\infty[t]$-module and $s_i$ for $N$ as a $\C_\infty[\sigma]$-module, then we have for $1\leq i\leq m$ and $1\leq j\leq \ell_i$
\[(t-\theta)^{j-1} e_i = s_{\ell_1 + \dots + \ell_i - j + 1}.\]
\end{enumerate}
If $N$ satisfies these two conditions, then for $\alpha = (\alpha_1,\dots,\alpha_r)^\top \in \TT_\theta^r$ we have
\begin{equation}\label{E:delta_0^N formula from NP21}
\delta_0^N(\alpha) =  
\begin{pmatrix}
\partial_t^{\ell_1-1}(\alpha_1)\\
\vdots\\
\partial_t(\alpha_1)\\
\alpha_1\\
\vdots\\
\partial_t^{\ell_r-1}(\alpha_r)\\
\vdots\\
\partial_t(\alpha_r)\\
\alpha_r\\
\end{pmatrix}\Bigg|_{t=\theta} .
\end{equation}

Similarly, by \cite[Remark 3.5.13]{NP21} we can give a similar construction for our $t$-motive $M$ so that the map $\delta_0^M$ can similarly be expressed using hyperderivatives.

We comment that the above constructions are commensurate with our Lemma \ref{L:Bases}, since if one changes the $\C_\infty[\sigma]$-basis for $N$ using a matrix $\alpha \in \Mat_{d}(\C_\infty[\sigma])$, then one can simply change the $\C_\infty[\tau]$-basis of $M$ by $(\alpha^*)\inv$ and our lemma still holds. A similar calculation also holds for a change of $\C_\infty[t]$-bases for $M$ and $N$.

We now state a sufficient condition for the matrices $Q_i$ and $P_i$ of Corollaries \ref{C:Exp coeff} and \ref{C:Log Coeff} to be invertible.

\begin{theorem}\label{T:Invertibility of P_i}
We maintain the conditions (1) and (2) given above for the dual $t$-motive $N$. We further require that the matrix $C$ be upper triangular. In this situation, the coefficients $P_i$ of the logarithm function of the associated $t$-module $\phi$ are invertible.
\end{theorem}

\begin{proof}
Our starting point is Corollary \ref{C:Log Coeff}. By definition \ref{D:Dual A Motive} we see that $\det(\Phi) = (t-\theta)^k$ for some $k$. Then, since $N$ is finitely generated as a $\C_\infty[\sigma]$-module, we may assume without loss of generality that each of the diagonal elements of $\Phi$ is a positive power of $(t-\theta)$. Thus
\[\Phi\inv = \begin{pmatrix}
\frac{1}{(t-\theta)^{k_1}} & * & \dots & *\\
0 & \frac{1}{(t-\theta)^{k_2}} & \dots & *\\
\vdots & \vdots & \ddots & \vdots \\
0 & 0 & \dots & \frac{1}{(t-\theta)^{k_r}}
\end{pmatrix},\]
for positive integers $k_m$. For each $1\leq f\leq d$, the basis element $h_f = (t-\theta)^k e_j$, for some non-negative $k$ and some $j$ with $1\leq j\leq r$. Viewing $N\isom \C_\infty[t]^r$, we may assume $e_j$ is the $j$th standard basis vector. Thus the $f$th column of $P_i$ is
\[\delta_0^N\left((\Phi\inv)\twist \cdots (\Phi \inv)\twistk{i}(t-\theta^{q^i})^{k}e_j\right ).\]
We note that the $j$th entry of $(\Phi\inv)\twist \cdots (\Phi \inv)\twistk{i}(t-\theta^{q^i})^{k}e_j$ equals
\begin{equation}\label{E:jth entry of Phi inv}
\frac{(t-\theta^{q^i})^{k}}{((t-\theta^q)\cdots (t-\theta^{q^i}))^{k_j}}.
\end{equation}
By \eqref{E:delta_0^N formula from NP21}, $\delta_0^N(\sigma^{-i}(h_f))$ is a block vector, whose blocks consist of successive $t$-hyperderivatives of functions, which are then evaluated at $t=\theta$. Combining this fact with our above formulas, we see that the $\ell_{1} + \dots + \ell_{j-1} + 1$ through $\ell_{1} + \dots + \ell_{j-1} + \ell_j$ entries of $\delta_0^N(\sigma^{-i}(h_f))$ are the $s=0$ through $s=\ell_j-1$ hyperderivatives of \eqref{E:jth entry of Phi inv} which are then evaluated at $t=\theta$. The final ingredients for our proof are the hyperderivative product formula \cite[Prop. 2.3.12]{Pap} and the formula given in \cite[Cor. 2.4.7]{Pap}, which states
\begin{equation}\label{E:derivative of t-theta}
\partial_t^j\left( \frac{1}{(t-\theta^{q^i})^m}\right ) = (-1)^j \binom{m+j-1}{j} \frac{1}{(t-\theta^{q^i})^{m+j}}.
\end{equation}
Taken together, these formulas allow us to compute the hyperderivatives of the terms appearing in \eqref{E:jth entry of Phi inv}. In conclusion, we see that $P_i$ is a block-upper-diagonal matrix, whose diagonal blocks all consist of successive hyperderivatives of terms of the form \eqref{E:jth entry of Phi inv}, which we can calculate using \eqref{E:derivative of t-theta}. Finally, a short calculation involving symmetric Pascal matrices (to account for the binomial coefficients appearing in \eqref{E:derivative of t-theta}) then shows that each of the blocks along the diagonal is invertible. We conclude that the matrix $P_i$ is invertible. Alternatively, we may conclude the proof by observing that the diagonal blocks of $P_i$ are each equal to the matrix coefficients the exponential function of the $\ell_j$th tensor power of the Carlitz module, which we can prove are invertible using formulas in \cite{Pap}.
\end{proof}

On the other hand, if our $t$-motive $M$ satisfies a condition similar to that given for $N$, then we can give a sufficient condition for the exponential coefficients $Q_i$ to be invertible. Namely:

\begin{enumerate}
\item There exists $C\in \text{GL}_r(\C_\infty[t])$ such that the $\tau$-action on $N$ is represented by a matrix
\[\Phi = C \begin{pmatrix}
(t-\theta)^{\ell_1} & & \\
 & \ddots &\\
  & & (t-\theta)^{\ell_r}
\end{pmatrix},\]
such that $\ell_k \geq 0$ with $\sum \ell_k = d$, the $\C_\infty[\tau]$-dimension of $N$.
\item If $e_i$ are basis elements for $N$ as a $\C_\infty[t]$-module and $s_i$ for $N$ as a $\C_\infty[\sigma]$-module, then we have for $1\leq i\leq m$ and $1\leq j\leq \ell_i$
\[(t-\theta)^{j-1} e_i = s_{\ell_1 + \dots + \ell_{i-1} + j}.\]
\end{enumerate}

\begin{theorem}\label{T:Invertibility of Q_i}
Let the $t$-motive $N$ satisfy conditions (1) and (2) given above. We further require that the matrix $C$ be upper triangular. In this situation, the coefficients $Q_i$ of the exponential function of the associated $t$-module $\phi$ are invertible.
\end{theorem}

\begin{proof}
The proof follows nearly identically to that of Theorem \ref{T:Invertibility of P_i}, using Corollary \ref{C:Exp coeff} rather than Corollary \ref{C:Log Coeff}.
\end{proof}

\begin{example}
We give a short example showing what can happen if our condition on the matrix $C$ in the above two theorems is not met. Let $N$ be a dual $t$-motive of $\C_\infty[t]$-rank 2 and $\C_\infty[\sigma]$-rank 3 with basis $\bh = \{h_1,h_2\}^\top$ such that
\[\sigma \bh = \begin{pmatrix}
0 & (t-\theta)^2 \\
(t-\theta) & 0
\end{pmatrix} \bh.\]
In this case $C = \begin{pmatrix}
0 & 1 \\
1 & 0
\end{pmatrix}$, which is not upper triangular. Using Anderson's functor produces a $3$-dimensional $t$-module $\phi$ with
\[\phi_t = \begin{pmatrix}
\theta & \tau & 0\\
0 & \theta & 1\\
\tau & 0 & \theta
\end{pmatrix}.\]
Our Corollary \ref{C:Log Coeff} then gives that the first couple logarithm coefficients are $P_0 = I$ (which is invertible), whereas
\[P_1 = \begin{pmatrix}
0 & \frac{1}{\theta - \theta^q} & \frac{1}{(\theta - \theta^q)^2}\\
\frac{-1}{(\theta - \theta^q)^2} & 0 & 0\\
\frac{1}{\theta - \theta^q} & 0 & 0\\
\end{pmatrix},\]
which is evidently not invertible. The pattern continues with the even coefficients being invertible and the odd coefficients being not invertible. We comment that our condition on $C$ is not necessary. For example, the exponential and logarithm coefficients of a large class of Drinfeld modules are invertible (being nonzero and dimension 1), but the associated $C$ matrices are not upper triangular.
\end{example}

\section{Product Formulas for Exp and Log}\label{S:Product Formulas}
In this section we prove formulas which express the exponential and logarithm functions as sums of infinite matrix products evaluated under $\delta_1^N$ and $\delta_{1}^M$. Pellarin proves somewhat similar looking formulas in the case of the Carlitz module in \cite[Thm. 1]{Pel12} and there are some hints towards the formulas of this section contained in the theorems of \cite[\S 6]{GP20} on trivial MZVs, but in general our formulas are completely original. As discussed in the introduction, we view a special case of these formulas (Corollary \ref{C:Carlitz zeta prod}) as giving a function field analogue of the Mellin transform formula for the Riemann zeta function. Thus, we are led to view the map $\delta_1^M$ as an analogue of integration. This analogy is explored in more detail in Remarks \ref{R:Mellin Remark} and \ref{R:Gamma Remark}.

\begin{definition}
For $n\in \Z_+$ we define a modified version of the $n$th partial sum of the pairings $F$ and $G$ from sections \S \ref{S:Exp Pairing} and \S \ref{S:Log Pairing}. For $x\in (L\otimes_{\F_q} A)[\tau]$ and $y \in (L\otimes_{\F_q} A)[\sigma]$ and $\bz \in L^d$ we define
\begin{equation}\label{D:Fpairing finite}
F_n(x,y;\bz) = \sum_{i=0}^n \sum_{k = 1}^{d} \delta_0^M(\tau^{-i}(x(g_k)))^\top \bz \sigma^i(y(h_k))) \in N
\end{equation}
\begin{equation}\label{D:Gpairing finite}
G_n(x,y) = \sum_{i=0}^n \sum_{k = 1}^{d} \delta_0^N\left(\sigma^{-i}(x(h_k))\right ) \tau^i(y(g_k))) \in M^d .
\end{equation}
\end{definition}

\begin{definition}
For a fixed $\bz \in L^d$, we extend the definition of $\delta_{1,\bz}^M$ to act coordinatewise on $M^d$. Then we have
\[F(x,y;\bz) = \lim_{n\to \infty} \delta_{1}^N(F_n(x,y;\bz)),\text{ and}\quad G(x,y;\bz) = \lim_{n\to \infty} \delta_{1,\bz}^M(G_n(x,y)),\]
where the second equality only holds for $\bz\in D\cap L^d$ as in Theorem \ref{T:Log from Pairing}.
\end{definition}

\begin{definition}\label{D:Theta tau}
For fixed $b\in A$, recall the definition of $\Theta_{b}$ from Lemma \ref{L:Bases} and write
\[\Theta_b = \Theta_0 + \Theta_1 \tau + \dots + \Theta_m \tau^m, \quad \Theta_i \in \Mat_d(L).\]
We then define
\begin{align*}
\Theta_{b,\tau^{m}} &=(\Theta_1 \tau + \dots + \Theta_m \tau^m) \\
\Theta_{b,\tau^{m-1}} &= (\Theta_2\twistinv \tau + \dots + \Theta_m\twistinv\tau^{m-1})\\
\vdots\\
\Theta_{b,\tau} &= (\Theta_m\twistk{1-m} \tau)
\end{align*}
and
\begin{align*}
\Theta_{b,\sigma^{m}} &=(\Theta_1 \tau + \dots + \Theta_m \tau^m) \\
\Theta_{b,\sigma^{m-1}} &= (\Theta_2 \tau + \dots + \Theta_m\tau^{m-1})\\
\vdots\\
\Theta_{b,\sigma} &= (\Theta_m \tau).
\end{align*}

\end{definition}

\begin{proposition}\label{P:Fn Gn linearity}
For the partial sums $F_n$ and $G_n$ we have
\begin{enumerate}
\item For all $a \in L$ we have
\[F_n(ax,y;\bz) = F_n(x,ay;\bz),\quad G_n(ax,y) = G_n(x,ay)\]
\item We have
\[F_n(x,\sigma y;\bz) - F_n(\tau x,y;\bz) = \sum_{k = 1}^{d} \delta_0^M(\tau^{-n}(x(g_k)))^\top \bz \sigma^{n+1}(y(h_k)))\]
\[G_n(x,\tau y) - G_n(\sigma x, y) =\sum_{k = 1}^{d} \delta_0^N\left(\sigma^{-n}(x(h_k)\right ) \tau^{n+1}(y(g_k)))\]
\item For $b\in A$, recall the definition of $\Theta_{b,\tau^{m-\ell}} \in \Mat_d(L[\tau])$ from Definition \ref{D:Theta tau} for $1\leq \ell\leq m$ and let $\be_k\in L^r$ be the $k$th standard basis vector. For all $x,y\in A$ we have
\[F_n(x,b y;\bz) - F_n(b x,y;\bz) = \sum_{k = 1}^{d}\sum_{\ell = 0}^{m-1} \delta_0^M(\tau^{\ell-n}(x(g_k)))^\top \bz \sigma^{n}(y(\be_k^\top\Theta_{b,\sigma^{m-\ell}}^*\bh))),\]
\[G_n(x,b y) - G_n(b x,y) = \sum_{k = 1}^{d}\sum_{\ell = 0}^{m-1} \delta_0^N(\sigma^{\ell-n}(x(h_k))) \tau^{n}(y(\be_k^\top\Theta_{b,\tau^{m-\ell}}\bg))),\]
\end{enumerate}
\end{proposition}

\begin{proof}
The proofs of parts (1) and (2) follow as in the proofs of Propositions \ref{P:F Bilinearity} and \ref{P:G Bilinearity}, and we leave these details to the reader. Part (3) follows similarly to the proof of Propositions \ref{P:F Bilinearity}(4) and \ref{P:G Bilinearity}(3), so we merely sketch the argument here. We begin by observing that
\[F_n(bx,y;\bz) =  \sum_{i=0}^n \sum_{k = 1}^{d} \delta_0^M(\tau^{-i}( x(\Theta_{b,k}\bg)))^\top \bz \sigma^i(y(h_k)) .\]
We then apply parts (1) and (2) of this proposition to the individual terms of $\Theta_{b,k}\bg$. Note that contrary to the proof of Proposition \ref{P:F Bilinearity}, $F_n$ is not $\tau$- and $\sigma$-symmetric; we get an extra term as observed in part (2) every time we transform a $\tau$ from the first coordinate into a $\sigma$ in the second coordinate. After a short calculation, we find that all the extra terms give
\[F_n(x,by;\bz) - F_n(bx,y;\bz) =  \sum_{k = 1}^{d}\sum_{\ell = 0}^{m-1} \delta_0^M(\tau^{\ell-n}(x(g_k)))^\top \bz \sigma^{n}(y(\be_k^\top\Theta_{b,\sigma^{m-\ell}}^*\bh))).\]
The proof of the formula for $G_n$ follows nearly identically.
\end{proof}

The next proposition gives a product formula for our pairings $F$ and $G$. We term this a product formula because powers of $\tau$ and $\sigma$ can be expressed as products of twists of matrices related to $\Phi$ and so as we take the limit as $n\to \infty$ (after proper normalization) we get a finite sum of infinite matrix products. We shall give explicit examples of such product formulas as they relate to special values of zeta and multiple zeta functions shortly.

\begin{theorem}[Product Formulas]\label{T:Product Formulas}
Let $b,x,y\in A$ and let $\Theta_{b,\tau}$ be as in Definition \ref{D:Theta tau}. For the pairings $F$ and $G$ from sections \S \ref{S:Exp Pairing} and \S \ref{S:Log Pairing} we have
\[F(x,y;\bz) = \lim_{n\to \infty} \delta_1^N\left((bI - d[b]^\top)\inv \sum_{k = 1}^{d}\sum_{\ell = 0}^{m-1} \delta_0^M(\tau^{\ell-n}(x(g_k)))^\top \bz \sigma^{n}(y(\be_k^\top\Theta_{b,\sigma^{m-\ell}}^*\bh)))\right ),\]
\[G(x,y;\bz) = \lim_{n\to \infty} \delta_{1,\bz}^M\left((bI - d[b])\inv \sum_{k = 1}^{d}\sum_{\ell = 0}^{m-1} \delta_0^N(\sigma^{\ell-n}(x(h_k))) \tau^{n}(y(\be_k^\top\Theta_{b,\tau^{m-\ell}}\bg)))\right ).\]
\end{theorem}

\begin{proof}
From the definition of $F_n(x,y;\bz)$ we see that for any $b\in A$
\[F_n(x,by;\bz) = b F_n(x,y;\bz).\]
Then, by the $A$-linearity of $\delta_0^M$ described in Proposition \ref{P:A linearity}, we find
\[F_n(bx,y;\bz) = d[b]^\top F_n(x,y;\bz).\]
Finally, by subtracting these, factoring and using Proposition \ref{P:Fn Gn linearity}(3) we deduce the first formula. The second formula follows identically.
\end{proof}

\begin{example}
In the case of the Carlitz module (for $d=1$), for example, we have
\begin{align*}
F(1,1;\bz) &= \lim_{n\to \infty} \delta_1^N\left (\frac{1}{t-\theta} \delta_0^M(\tau^{-n}(1))\cdot z \cdot \sigma^{n+1}(1)\right )\\
&= \lim_{n\to \infty} \delta_1^N\left (z\frac{(t-\theta^{1/q})\cdots (t-\theta^{1/q^n})}{(\theta-\theta^{1/q})\cdots (\theta - \theta^{1/q^n})}\right )
\end{align*}
\end{example}

\begin{corollary}\label{C:Exp and Log Products}
Let $\phi$ be a $\Xi$-regular Anderson $A$-module.
\begin{enumerate}
\item For all $\bz \in \C_\infty^d$ we have
\[\Exp_\phi(\bz) = \lim_{n\to \infty} \delta_1^N\left((bI - d[b]^\top)\inv \sum_{k = 1}^{d}\sum_{\ell = 0}^{m-1} \delta_0^M(\tau^{\ell-n}(x(g_k)))^\top \bz \sigma^{n}(y(\be_k^\top\Theta_{b,\sigma^{m-\ell}}^*\bh)))\right ),.\]
\item For all $\bz$ in the domain of convergence of $\Log_\phi$ we have
\[\Log_\phi(\bz) = \lim_{n\to \infty} \delta_{1,\bz}^M\left((bI - d[b])\inv \sum_{k = 1}^{d}\sum_{\ell = 0}^{m-1} \delta_0^N(\sigma^{\ell-n}(x(h_k))) \tau^{n}(y(\be_k^\top\Theta_{b,\tau^{m-\ell}}\bg)))\right ).\]
\end{enumerate}

\end{corollary}
\begin{proof}
This follows immediately from Theorems \ref{T:Exponential from Pairing}, \ref{T:Log from Pairing} and \ref{T:Product Formulas}.
\end{proof}

\begin{example}\label{Ex: Tensor powers Carlitz Prod Formula}
We return to the notation of Example \ref{E:Tensor Power of Carlitz} for the $n$th tensor power of the Carlitz module. In this situation, the formula given in Theorem \ref{T:Product Formulas} is especially clean and interesting. In this setting we have
\begin{align}\label{E:Gm(1,1) formula}
G_m(1,1) &= \sum_{i=0}^m \sum_{k=1}^n \delta_0^N(\sigma^{-i}((t-\theta)^{n-k}))\tau^i((t-\theta)^{k-1})\nonumber \\
&= \sum_{i=0}^m \sum_{k=1}^n \delta_0^N\left(\frac{(t-\theta^{q^i})^{n-k}}{(t-\theta^{q})^n \cdots (t-\theta^{q^i})^n}\right )(t-\theta)^n\cdots (t-\theta^{q^{i-1}})^n(t-\theta^{q^i})^{k-1} 
\end{align}
From equation \eqref{E:Carlitz Def} we see that $\Theta_{t,\tau}$ on the basis $\{g_1,\dots,g_n\}$ is given by
\[\Theta_{t,\tau} = \begin{pmatrix}
0 & \cdots & 0\\
\vdots & & \vdots \\
\tau & \cdots & 0
\end{pmatrix}.\]
From this, we see that most of the terms in Theorem \ref{T:Product Formulas} vanish and we are left with
\begin{align}\label{E:G11 formula}
G(1,1;\bz) &=  \delta_{1,\bz}^M\left(\lim_{m\to \infty}(t - d[t])\inv \delta_0^N(\sigma^{-m}(1)) \tau^{m+1}(1)\right )\\
\nonumber & = \delta_{1,\bz}^M\left(\lim_{m\to \infty} (t - d[t])\inv \delta_0^N\left(\frac{1}{(t-\theta^{q})^n \cdots (t-\theta^{q^m})^n}\right ) (t-\theta)^n\cdots (t-\theta^{q^{m}})^n\right ).
\end{align}
\begin{remark}
We comment that the map $\delta_{1,\bz}^M$ in the above formula is still the extension described in Definition \ref{D:N[[sigma]]}; in the second equality one must implicitly use the isomorphism $M \isom \C_\infty[t] \isom \C_\infty[\tau]^d$.
\end{remark}
Recall the formula (see \cite[2.2.1]{BP20}) for the fundamental period of the Carlitz module $\tilde \pi$
\begin{equation}\label{E:Carlitz period def}
\tilde \pi = \theta(-\theta)^{1/(q-1)} \prod_{i=1}^\infty \Bigl(
1 - \frac{\theta}{\theta^{q^i}} \Bigr)^{-1} \in \C_\infty
\end{equation}
and the definition of the Anderson-Thakur function 
\begin{equation}\label{E:And-Thak Function}
\omega_C = (-\theta)^{1/(q-1)} \prod_{i=0}^\infty \biggl( 1 - \frac{t}{\theta^{q^i}} \biggr)^{-1} \in \TT
\end{equation}
from \cite{AT90}. Finally, recall that the bottom coordinate of the map $\delta_0^N$ in this situation is merely evaluation at $\theta$ and that $(t-d[t])$ is upper triangular with bottom-right coordinate $(t-\theta)$. Putting this all together, and examining the bottom coordinate of equation \eqref{E:G11 formula} we conclude the following theorem.
\end{example}

\begin{theorem}\label{T:Carlitz Prod Formula}
Let $\Log_C\on$ denote the logarithm function for the $n$th tensor power of the Carlitz module, let $M$ be the corresponding $t$-motive and let $\delta_{1,\bz}^M$ be the map from Definition \ref{D:delta^M defs} for any $\bz \in \C_\infty^n$ inside the domain of convergence for $\Log_C\on$. Then, if we let $p_n$ denote the projection onto the $n$th coordinate, we have the formula
\[p_n(\Log_C\on(\bz)) = \delta_{1,\bz}^M\left((-1)^{n-1}\frac{\tilde \pi^n}{(\theta-t)\omega_C^n}\right ).\]
\end{theorem}

We denote the Carlitz zeta value by
\begin{equation}\label{E:zeta def}
\zeta_A(n) = \sum_{\substack{a\in \F_q[\theta] \\ a\text{ monic}}} \frac{1}{a^n}, \quad n\in \N.
\end{equation}
Denote $D_{0}:=1$ and $D_{i}:=\prod_{j=0}^{i-1}(\theta^{q^{i}}-\theta^{q^{j}})\in \C_\infty$ for $i\in \N$. For any non-negative integer $n$, recall that the Carlitz gamma value is defined by
\begin{equation}\label{E: Carlitz factorials}
\Gamma(n+1):=\prod_{i=0}^{\infty}D_{i}^{n_i}\in A 
\end{equation}
where the $n_i \in \Z_{\geq 0}$ are given by writing the base $q$-expansion $n=\sum_{i=0}^{\infty} n_{i}q^{i}$ for $0\leq n_i\leq q-1$.

We now make an application of Theorem \ref{T:Carlitz Prod Formula} to Carlitz zeta values using the formulas for $\zeta_A(n)$ from \cite[Thm. 3.8.3]{AT90}. Let $\bz  \in K^n$ be the special point of \cite[3.8.2]{AT90} (in our notation $\bz = \delta_1^N(H_n(t))$ where $H_n(t)$ is the $n$th Anderson-Thakur polynomial; see \cite[\S 4.8]{GND20b}) and suppose that $n$ and $q$ are such that $\bz$ is within the domain of convergence of $\Log_{C\on}$. Further, we apply \cite[Thm. 6.7]{Gre17a} to the case of Anderson generating functions (see \cite{EGP14} for details) for the $n$th tensor power of the Carlitz module to find
\begin{equation}\label{E:And gen for Catlitz tensor powers}
\Exp_{C\on}\left ( (d[t]-tI)\inv \Pi_n\right ) = ((t-\theta)^{n-1}\omega_C^n,\dots,(t-\theta)\omega_C^n, \omega_C^n)^\top,
\end{equation}
where $\Pi_n\in \C_\infty^n$ is a fundamental period for $\Exp_{C\on}$. We note that in the above formula, we are viewing $t$ as a central variable which is unaffected by Frobenius twisting, in the sense of the formula for $\omega(t)$ in \cite[\S 2.1]{AP15} Let us denote $\bu = (d[t]-tI)\inv \Pi_n$, so that we have
\begin{equation}\label{E:Last coord of And Gen}
p_n(\Exp_{C\on}(\bu)) = \omega_C^n,
\end{equation}
and also let $u = \frac{\tilde \pi^n}{\theta-t}$ (by \cite[Thm. 6.7]{Gre17a} $u$ is the bottom coordinate of $\bu$).

\begin{corollary}\label{C:Carlitz zeta prod}
Let $\bz$ be the special point of \cite[3.8.2]{AT90} subject to the notation and restrictions given above. Then we have
\[\delta_{1,\bz}^M\left((-1)^{n-1}\frac{\tilde \pi^n}{(\theta-t)\omega_C^n}\right ) = \delta_{1,\bz}^M\left((-1)^{n-1}\frac{u}{p_n(\Exp_{C\on}(\bu))}\right ) = \Gamma(n)\zeta_A(n) ,\]
where $M$ is the $t$-motive corresponding to $C\on$.
\end{corollary}

\begin{proof}
The first equality follows from the discussion above. The second equality follows immediately from Theorem \ref{T:Carlitz Prod Formula} and by Theorem 3.8.3 of \cite{AT90}.
\end{proof}

\begin{remark}
It is possible to give Corollary \ref{C:Carlitz zeta prod} for all values of $n\geq 1$ without constraints on the size of $\bz$. However, in general the statement becomes much more cumbersome due to the fact that $\bz$ is no longer in the domain of convergence of the logarithm. One would then need to break the LHS of Corollary \ref{C:Carlitz zeta prod} into a $C\on$-linear combination of logarithms, as is done in \cite[3.9]{AT90}, and the left hand side of Corollary \ref{C:Carlitz zeta prod} would then become a sum of terms of the form $\delta_{1,C_\alpha\on \bz}^M(\cdots)$. We comment that the simple version of Corollary \ref{C:Carlitz zeta prod} holds at least for $1\leq n \leq q$. In general, the exact values for which Anderson and Thakur's special point $\bz$ lies in the domain of convergence of the logarithm is a subtle question which we leave for future study.
\end{remark}

\begin{remark}\label{R:Mellin Remark}
As discussed in the introduction, we view Corollary \ref{C:Carlitz zeta prod} as a function field analogue of the Mellin transform formula for the Riemann zeta function
\[\int_0^\infty \frac{1}{e^x-1}x^{s-1}dx = \Gamma(s)\zeta(s),\]
where the map $\delta_{1,\bz}^M$ serves as a replacement of integration. It is likely that the more natural way to state Theorem \ref{T:Carlitz Prod Formula} and Corollary \ref{C:Carlitz zeta prod} is using $n$-dimensional vectors. However, we have chosen to state the theorem after projecting onto the last coordinate to highlight the analogy between our formulas and the Mellin transform formula.
\end{remark}

\begin{remark}\label{R:Gamma Remark}
Here we briefly describe how we can recover function field gamma values using $\delta_{1,\bz}^M$. Set $\bz = \delta_1^N(H_n)$ ($H_n$ is an Anderson-Thakur polynomial \cite[\S 4.8]{GND20b}) which is Anderson and Thakur's special point. Then a short calculation gives
\[p_n(\delta_{1,\bz}^M(\delta_0^N(1)\cdot (t-\theta)^{n-1})) = \Gamma(n).\]
Note that this is the first term of the expansion given in \eqref{E:Gm(1,1) formula}.
\end{remark}

\begin{example}\label{Ex:MZV (3,1)}
In this example we sketch an application of our theorems to a specific $t$-module whose logarithm evaluates to give Carlitz multiple zeta values (MZVs). There are substantial details which go into this construction which are explained elsewhere, thus here we merely provide a sketch and describe how our theorems relate to it (see particularly \cite[Ex. 5.4.2]{CM19b} and \cite[Ex. 6.2]{GND20b}). We emphasize that this is a specific example of a general phenomenon; it should be possible to use this approach for a general MZV and we plan to explore this theme in more detail in a future project.

For any $r$-tuple of positive integers $\fs=(s_{1},\ldots,s_{r})\in \N^{r}$, define the MZV for the tuple $\fs$ by
\begin{equation}\label{D:MZVdef}
\zeta_{A}(\fs):=\sum \frac{1}{a_{1}^{s_1}\cdots a_{r}^{s_r}}\in K_{\infty},
\end{equation}
where the sum is over all $r$-tuples of monic polynomials $a_{1},\ldots,a_{r}$ in $A$ with the restriction $|a_{1}|_{\infty}> |a_{2}|_{\infty}>\cdots>|a_{r}|_{\infty}$.

We follow the construction in \cite[\S 5.1]{GND20b} of star dual $t$-motives. We maintain the notation of the previous example for $A=\F_q[t]$, $K=\F_q(t)$, an so on and fix an $r$-tuple $\fs = (s_1,\dots,s_r)\in \Z_+^r$. We define a dual $t$-motive with $\sigma$-action given on a $\C_\infty[t]$-basis $\{h_1,\dots,h_r\}$ by 
\begin{align*}
\Phi^\star &:=
\begin{pmatrix}
\Phi^\star_{1,1} & &  \\
\vdots & \ddots & \\
\Phi^\star_{r,1} & \dots & \Phi^\star_{r,r}
\end{pmatrix} \in \Mat_{r}(\overline K[t])
\end{align*}
where for $1 \leq \ell \leq j \leq r$, 
\begin{equation} \label{eq: Phi star}
\Phi^\star_{j,\ell}=(-1)^{j-\ell} \prod_{\ell \leq k <j} H_k^{(-1)} (t-\theta)^{s_\ell+\dots+s_r},
\end{equation}
and where $H_k\in K[t]$ is a particular Anderson-Thakur polynomial coming from the $r$-tuple $\fs$ (see \cite[\S 4.8]{GND20b}). As a particular example, if we set $q=2$ and set $\fs = (1,3)$, then
\[\Phi^\star = \begin{pmatrix}
(t-\theta)^4 & 0\\
-(t-\theta)^5 & (t-\theta)
\end{pmatrix}\]
and the associated $t$-module is then given by
\begin{align*}
\phi_t=\left(
\begin{array}{c c c c | c  }
\theta & 1 &  &  &  \\
& \theta & 1 &  &  \\
&  & \theta & 1 &  \tau \\
\tau &  &  & \theta & (\theta^2+\theta)\tau \\
\hline
&  &  &  & \theta +\tau
\end{array}
\right).
\end{align*}
Let $e_i$ for $1\leq i\leq 5$ be the standard basis vectors, so that $\{e_i\}$ is a $\C_\infty[\tau]$- and $\C_\infty[\sigma]$-module basis for $M$ and $N$, respectively. A short calculation shows that $m_1=e_1$ and $m_2 = e_5$ forms a $\C_\infty[t]$-basis for $M$ and $n_1=e_4$ and $n_2 = e_5$ forms a $\C_\infty[t]$-basis for $N$. A quick check shows that these bases satisfy the conditions of Lemma \ref{L:Bases}. Further, we calculate that we can express 
\begin{align*}
e_1 = m_1, \quad & e_1 = (t-\theta)^3 n_1\\
e_2 = (t-\theta)m_1,  \quad & e_2 = (t-\theta)^2 n_1\\
e_3 = (t-\theta)^2 m_1,  \quad & e_3 = (t-\theta) n_1\\
e_4 = (t-\theta)^3 m_1 + (t-\theta)m_2,  \quad & e_4 = n_1\\
e_5 = m_2,  \quad & e_5 =  n_2.
\end{align*}
So, the first few terms $G(1,1)$ are given by
\[G(1,1) = \delta_1^M(\delta_0^N((t-\theta)^3 n_1)m_1  + \delta_0^N((t-\theta)^2 n_1) (t-\theta)m_1 + \delta_0^N( (t-\theta)n_1)(t-\theta)^2m_1+ \cdots\]
In order to apply Theorem \ref{T:Product Formulas} we calculate that $\Theta_t = \phi_t$, and hence
\[\Theta_{t,\tau} = \left(
\begin{array}{c c c c | c  }
 &  &  &  &  \\
&  &  &  &  \\
&  &  &  &  \tau \\
\tau &  &  &  & (\theta^2+\theta)\tau \\
\hline
&  &   &  & \tau
\end{array}
\right).\]
Thus, Theorem \ref{T:Product Formulas} gives
\begin{align*}
G(1,1) &= \delta_{1}^M\bigg(\lim_{i\to \infty} (t-d[t])\inv \big[\delta_0^N(\sigma^{-i}((t-\theta)n_1)) \tau^{i+1}(m_2) \\ &+\delta_0^N(\sigma^{-i}(n_1))[\tau^{i+1}(m_1) + (\theta^2 + \theta)^{q^i}\tau^{i+1}(m_2)] + \delta_0^N(\sigma^{-i}(n_2))\tau^{i+1}(m_2)\big ]\bigg ).
\end{align*}
Collecting like terms and writing $m_1,n_1 = (1,0)^\top$ and $m_2,n_2 = (0,1)^\top$ gives
\[G(1,1) = \delta_{1}^M\bigg(\lim_{i\to \infty} (t-d[t])\inv \bigg[\delta_0^N\left(\sigma^{-i}\begin{pmatrix}t-\theta^2\\ 1 \end{pmatrix}\right ) \tau^{i+1}\begin{pmatrix}
0\\1
\end{pmatrix} + \delta_0^N\left(\sigma^{-i}\begin{pmatrix}
1\\0
\end{pmatrix}\right )\tau^{i+1}\begin{pmatrix}
1\\0
\end{pmatrix}\bigg ]\bigg ).\]
Finally, if we set $\bz = (0,0,0,0,-1)^\top$ and if we let $p_4$ denote the projection on to the $4$th coordinate, then by \cite[Ex. 6.2]{GND20b} we find
\begin{equation}\label{E:MZV product formula}
p_4(G(1,1,\bz)) = (\theta^2+\theta) \zeta_A(1,3),
\end{equation}
thus giving us a ``product formula" for $\zeta_A(1,3)$. The above strategy applies equally well to $\zeta_A(\fs)$ for any $r$-tuple $\fs$, but it becomes increasingly difficult to calculate the individual terms involved in the left hand side of \eqref{E:MZV product formula}. Studying this phenomenon for a general $r$-tuple will be the subject of a future project.
\end{example}

\begin{remark}\label{R:MZV remark}
We observe that certain zeta values seem to appear in the terms of our product formula above. For example, the term
\[\delta_0^N\left(\sigma^{-i}\begin{pmatrix}
1\\0
\end{pmatrix}\right )\tau^{i+1}\begin{pmatrix}
1\\0
\end{pmatrix} = \delta_0^N\left((\Phi^{\star\top(-1)}\cdots \Phi^{\star\top(-i)})\inv\begin{pmatrix}
1\\0
\end{pmatrix}\right )(\Phi^{\star}\cdots \Phi^{\star(i)})\begin{pmatrix}
1\\0
\end{pmatrix}.\]
Examining the formula for $\Phi^\star$ given in Example \ref{Ex:MZV (3,1)}, we see that the 4th coordinate of the above expression equals
\[\frac{1}{(\theta-\theta^{q})^4 \cdots (\theta-\theta^{q^i})^4}
(t-\theta)^4\cdots (t-\theta^{q^{i-1}})^4.
\]
This matches the formula we obtained in Example \ref{Ex: Tensor powers Carlitz Prod Formula} for $n=4$, which hints that a term consisting of $\zeta_A(4)$ appears in the expression for $p_4(G(1,1,\bz))$ (the $\delta_1$ maps are different in each case, so we can't make a direct comparison). We are hopeful that studying formulas of this type could lead to new expressions for $\F_q(\theta)$-linear relations between multiple zeta values. However, the $*$ coordinate complicates the situation. This theme will be explored in a future project.
%\[p_4 \left(\delta_{1}^M\bigg(\lim_{i\to \infty} (t-d[t])\inv \delta_0^N\left(\sigma^{-i}\begin{pmatrix}
%1\\0
%\end{pmatrix}\right )\tau^{i+1}\begin{pmatrix}
%1\\0
%\end{pmatrix} \bigg)\right )
%=
%\begin{pmatrix}
%\Gamma_4 \zeta_A(4)\\
%*
%\end{pmatrix}
%\]
\end{remark}

\section{Log-Algebraic Criterion}\label{S:Log Alg}
In this section we make a short application of the motivic pairings of Sections \ref{S:Exp Pairing} and \ref{S:Log Pairing} to give a new log-algebraicity criterion. The main theorem of this section is Theorem \ref{T:Main Log Alg Theorem}, which gives a generalization of an unpublished theorem originally due to Anderson (see \cite[Cor. 2.5.23]{HJ20}). Although the main result of this section merely recovers (and slightly extends) a known theorem, the methods of proof we use here are new. The pairings $H$ and $I$ which we define in this section appear in a concurrent project which develops a non-commutative factorization of the exponential function, and thus we are inspired to include them here.

We now develop two new pairings, $H$ and $I$, which relate the $F$ and $G$ pairings. The $I$ pairing will prove especially crucial in this section. Throughout this section, we set $L= \C_\infty$.

\begin{definition}
Fix an integer $\ell\geq 0$, and let $x\in (L\otimes A)[\tau]$ and $y\in (L\otimes A)[\sigma]$. Recall that $\{g_i\}$ and $\{h_i\}$ are the bases of Lemma \ref{L:Bases}. Then define the pairing
\[H_\ell:(L\otimes A)[\tau]\times (L\otimes A)[\sigma]\to \Mat_{d\times d}(L),\]
by setting
\[H_\ell(x,y) = \sum_{j=0}^\ell\sum_{k=1}^d \delta_0^N(\sigma^{-j}(y(h_k))) \left(\delta_0^M(\tau^{j-\ell}(x(g_k)))^\top\right )\twistk{\ell}.\]
\end{definition}

We also define a finer version of the pairing.

\begin{definition}
For, $\ell\geq 0$, for $x\in (L\otimes A)[\tau]$ and $y\in (L\otimes A)[\sigma]$, and for any $g\in M$ and $h\in N$, we define
\[H_\ell(x,y;g,h) = \sum_{j=0}^\ell \delta_0^N(\sigma^{-j}(y(h))) \left(\delta_0^M(\tau^{j-\ell}(x(g)))^\top\right )\twistk{\ell}.\]
\end{definition}

\begin{lemma}\label{L:H fine linearity}
The pairing $H_\ell(x,y;g,h)$ satisfies
\begin{enumerate}
\item $H_\ell(\tau x,y;g,h) = H_\ell(x,\sigma y;g,h)$
\item $H_\ell(a x,y;g,h) = H_\ell(x,ay;g,h)$ for all $a\in L$
\end{enumerate}
\end{lemma}

\begin{proof}
The proof for (1) follows by recalling that $\delta_0(\tau(g)) = 0 = \delta_0(\sigma(h))$ for all $g\in M$ and $h\in N$. The proof for (2) follows using the $L$-commutativity relations for $\tau$ and $\sigma$ and the fact that both maps $\delta_0^M, \delta_0^N$ are $L$-linear.
\end{proof}

\begin{proposition}\label{P:H pairing bilinearity}
The pairing $H_\ell(x,y)$ satisfies
\begin{enumerate}
\item $H_\ell(\tau x,y) = H_\ell(x,\sigma y;g,h)$
\item $H_\ell(a x,y) = H_\ell(x,ay)$ for all $a\in L$
\item $H_\ell(b x,y) = H_\ell(x,by)$ for all $b\in A$ and for $x,y\in A$
\end{enumerate}
In particular, $H_\ell(x,y) = 0$ whenever $\ell>0$ and $x,y\in A$ and $H_0(1,1) = I_d$.
\end{proposition}

\begin{proof}
Parts (1) and (2) follow directly from Lemma \ref{L:H fine linearity}. The proof of part (3), follows similarly to the proof of Proposition \ref{P:F Bilinearity} part (4); we leave the details to the reader. To prove the last statement of the proposition, we note that for $b\in A$ we have
\begin{align*}
H_\ell(x,by) & = \sum_{j=0}^\ell\sum_{k=1}^d \delta_0^N(\sigma^{-j}(by(g_k))) \left(\delta_0^M(\tau^{j-\ell}(c(h_k)))^\top\right )\twistk{\ell}\\
&= \sum_{j=0}^\ell\sum_{k=1}^d d[b]\delta_0^N(\sigma^{-j}(y(g_k))) \left(\delta_0^M(\tau^{j-\ell}(x(h_k)))^\top\right )\twistk{\ell}\\
&= d[b]H_\ell(x,y)\\
\end{align*}
whereas
\begin{align*}
H_\ell(bx,y) & = \sum_{j=0}^\ell\sum_{k=1}^d \delta_0^N(\sigma^{-j}(y(g_k))) \left(\delta_0^M(\tau^{j-\ell}(bx(h_k)))^\top\right )\twistk{\ell}\\
&= \sum_{j=0}^\ell\sum_{k=1}^d \delta_0^N(\sigma^{-j}(y(g_k))) \left((d[b]^\top\delta_0^M(\tau^{j-\ell}(x(h_k))))^\top\right )\twistk{\ell}\\
&= \sum_{j=0}^\ell\sum_{k=1}^d \delta_0^N(\sigma^{-j}(y(g_k))) \left(\delta_0^M(\tau^{j-\ell}(x(h_k)))^\top\right) d[b])\twistk{\ell}\\
&= H_\ell(x,y) d[b]\twistk{\ell}.\\
\end{align*}
By part (3), $H_\ell(bx,y) = H_\ell(x,by)$ for $x,y\in A$, then a short linear algebra argument shows that 
\[d[b]H_\ell(x,y) =  H_\ell(x,y) d[b]\twistk{\ell}\]
can never happen if $\ell \geq 1$, unless $H_\ell(x,y)=0$. On the other hand, if $x=y=1$ and $\ell=0$, then $H_0(1,1) = I_d$.
\end{proof}

We also define a pairing which reverses the role of $M$ and $N$ from the previous pairing $H_\ell$. We define
\[I_\ell:(L\otimes A)[\tau]\times (L\otimes A)[\sigma]\to L,\]
by setting
\[I_\ell(x,y) = \sum_{j=0}^\ell\sum_{k,m=1}^d  \left(\delta_0^M(\tau^{-j}(x(g_k)))^\top\right )\twistk{j}\delta_0^N(\sigma^{j-\ell}(y(h_m)))\twistk{j}.\]

We also define a finer version of this pairing for $g \in M$ and $h\in N$ by setting

\[I_\ell(x,y,g,h) = \sum_{j=0}^\ell  \left(\delta_0^M(\tau^{-j}(x(g)))^\top\right )\twistk{j}\delta_0^N(\sigma^{j-\ell}(y(h)))\twistk{j}.\]

\begin{proposition}
For $I_\ell(x,y,g,h)$ (and consequently also for $I_\ell(x,y)$) we have
\begin{enumerate}
\item $I_\ell(bx,y) = I_\ell(x,by)$ for all $b\in A$
\item $I_\ell(\tau x,y) = I_{\ell-1}(x,y)\twist$
\item $ I_\ell(x,\sigma y) = I_{\ell-1}( x,y)$
\item $I_\ell(x,ay) = I_\ell(a^{q^\ell}x,y)$ for all $a\in L$
\end{enumerate}
\end{proposition}

\begin{proof}
The proof follows very similarly to the proof of Proposition \ref{P:H pairing bilinearity}. We leave the details to the reader.
\end{proof}

\begin{proposition}\label{P:I vanishing}
We have $I_\ell(1,1,g_k,h_m)=1$ when $k=m$ and $\ell=0$ and we have $I_\ell(1,1,g_k,h_m)=0$ otherwise.
\end{proposition}

\begin{proof}
Let $\bz = (z_1,\dots,z_d)^\top$ be in the domain of convergence of $\Log_\phi$. We first observe that by Theorems \ref{T:Exponential from Pairing} and \ref{T:Log from Pairing} we have
\begin{align*}
\Exp_\phi(\Log_\phi(\bz)) &= \delta_1^N\left(\sum_{i=0}^\infty \sum_{k = 1}^{d} \delta_0^M(\tau^{-i}(g_k))^\top \left[ \Log_\phi(\bz)\right ] \sigma^i(h_k))  \right )\\
&= \delta_1^N\left(\sum_{i=0}^\infty \sum_{k = 1}^{d} \delta_0^M(\tau^{-i}(g_k))^\top \left[ \delta_0^N\left(\sum_{j=0}^\infty \sum_{m = 1}^{d} \sigma^{-j}(h_m) \delta_{1,\bz}^M(\tau^j(g_m))  \right )\right ] \sigma^i(h_k))  \right )\\
&= \sum_{i=0}^\infty \sum_{k = 1}^{d}\sum_{j=0}^\infty \sum_{m = 1}^{d} (\delta_0^M(\tau^{-i}(g_k))\twisti)^\top   \delta_0^N\left(\sigma^{-j}(h_m)   \right )\twisti \delta_{1,\bz}^M(\tau^j(g_m))^{q^i}\delta_1^N\left(\sigma^i(h_k))  \right )\\
&= \sum_{i,j=0}^\infty \sum_{k,m = 1}^{d} (\delta_0^M(\tau^{-i}(g_k))\twisti)^\top   \delta_0^N\left(\sigma^{-j}(h_m)   \right )\twisti z_m^{q^{i+j}}\be_k 
\end{align*}
where $\be_k$ is the $k$th standard basis vector. Then, since $\Exp_\phi(\Log_\phi(\bz)) = \bz$, we compare this coordinate-wise with the above equality to conclude the proof.
\end{proof}

Because $(\C_\infty \otimes_{\F_q} A)$ is dense in $\YY_\theta$, it follows that every $h\in N_\Gamma$ (see Def. \ref{D:Tate algebras} and following discussion) we can write $\sigma^w(h)$ as a unique infinite $\C_\infty$-linear combination of the $\C_\infty[\sigma]$-basis elements of Lemma \ref{L:Bases},
\begin{equation}\label{E:sigma^i expression}
\sigma^w(h)  = c_{1,0}^wh_1 + c^w_{2,0}(h_2) + \dots  + c^w_{d,w}\sigma^w(h_d) + \mathcal O(\sigma^{w+1}),
\end{equation}
for $c_{i,j}^w \in \C_\infty$ ($w$ is simply an additional index for the coefficient, not a power) where $\mathcal O(\sigma^k)$ denotes terms with at least a $k$th power of sigma. Note in particular that if $h\in N$, then $c_{i,j}^w=0$ for $j<w$. Then, for each $w \geq 0$, we obtain an expression for $h$ by taking $\sigma^{-w}$ of both sides of \eqref{E:sigma^i expression} using Definition \ref{D:tau inverse}, which we label as
\begin{align}\label{E:h ell expression}
\begin{split}
h^w  =& d_{1,-w}^w \sigma^{-w}(h_1) + d_{2,-w}^w\sigma^{-w}(h_2) + \dots + d_{1,0}^w(h_1)+ \dots + d_{d,0}^w(h_d) + \mathcal O(\sigma).
\end{split}
\end{align}
We comment that if the coordinates of $h$ are all regular at $\Xi\twisti$ for $i\geq 1$, then each $d_{i,j}=0$ for $j<0$, since the coordinates of $(\Phi\inv)\twisti$ necessarily have poles at $\Xi\twisti$. Therefore, if $h \in N \subset N_\Gamma$, then the expressions $h^w$ are all identical for $w\geq 0$. Finally, for such an expression, we denote
\[\bd_w := 
\left (\begin{matrix}
d_{1,0}^w\\
d_{2,0}^w\\
\vdots\\
d_{n,0}^w
\end{matrix}\right ).\]

\begin{definition}\label{D:delta_* def}
For a function $h\in N_\Gamma$, if $\lim_{w\to \infty} \bd_w$ exists, then we denote this limit by $L_h$. Then let $N_\Gamma^* \subset N_\Gamma$ be the set of functions $h$ such that $L_h$ exists and is finite. We define a new extension of $\delta_0^N$ to $N_\Gamma^*$, which we call $\delta_*^N:N_\Gamma^* \to \C_\infty^d$, by setting
\[\delta_*^N(h) = \lim_{w\to \infty} \bd_w.\]
Note that if each coordinate of $h$ is regular at $\Xi\twisti$ for $i\geq 1$, then $\delta_*^N(h) = \delta_0^N(h)$.
\end{definition}

\begin{proposition}
The set $N_\Gamma^*$ is nonempty. In particular we have
\begin{enumerate}
\item $N \subset N_\Gamma^*$
\item If $h\in N_\Gamma$ and satisfies $h -\sigma(h)= \sigma(g)$, for some $g\in N$, then $h\in N_\Gamma^*$.
\end{enumerate}
\end{proposition}

\begin{proof}
By the discussion above, if $h\in N$, then $\delta_*^N(h) = \delta_0(h)$ and the limit in Definition \ref{D:delta_* def} exists. If $h\in N_\Gamma$ satisfies $h -\sigma(h)= \sigma(g)$, then we may write $\sigma^i(h) = h - \sigma(g) - \sigma^2(g) - \dots - \sigma^i(g)$. This allows us to calculate expression \eqref{E:sigma^i expression} for $h$ in terms of $\sigma^{-i}(h)$ and a short calculation shows that the limit in Definition \ref{D:delta_* def} exists and that $\delta_*^N(h) = \delta_1^N(g)$.
\end{proof}

\begin{definition}
For $1\leq k\leq d$, we define a function $J_k:N_\Gamma \to L$ similar to our pairings $F,G,H$ and $I$ of the previous sections. For $h \in N_\Gamma$ and $1\leq k\leq d$ we set
\[J_k(h) = \sum_{i=0}^\infty \left(\delta_0^M(\tau^{-i}(g_k))^\top\right )\twisti\delta_0^N(\sigma^{i}(h))\twisti\in L .\]
We then define $J:N_\Gamma \to L^d$ by
\[J(h) = (J_1,\dots,J_d)^\top.\]
\end{definition}

\begin{lemma}\label{L:h convergent}
For $h\in N_\Gamma$, the series $J(h)$ converges in $\C_\infty^d$.
\end{lemma}

\begin{proof}
%As explained in the proof of Proposition \ref{P:delta_0 extensions}, the map $\delta_0^N$ factors through $N/J^dN$. 
We temporarily view $N$ as a free $\C_\infty[t]$-module of rank equal to $r\cdot \ell$, where $\ell$ is the degree of $K/\F_q(t)$. The $t$-motive $N$ produces a $t$-module $\phi$ (which is equal to the Anderson $A$-module $\phi$ restricted to $\F_q[t]$ and has the same exponential function), thus we may assume without loss of generality that we are in the situation outlined in conditions (1) and (2) of Theorem \ref{T:Invertibility of P_i}, where the $\delta_0^N$ map is given by hyperderivatives evaluated at $\theta$. Then, a short calculation shows for $h \in N_\Gamma$, that $|\delta_0^N(h) |\leq |h|_\theta$. It follows that
\[|\delta_0^N(\sigma^i(h))| \leq |\sigma^i(h)|_\theta = |\Phi \cdots \Phi \twistk{1-i} h\twistk{-i}|_\theta \leq |\Phi|^{(1-1/q^i)/(1-1/q)}_\theta |h|^{1/q^i}_\theta.\]
The last quantity in the above inequality is bounded independently from $i$, and thus the convergence of the series $J(h)$ follows from the fact that the exponential function is entire and from the fact that $\delta_0^M(\tau^{-i}(g_k))$ appears in the coefficients of  $\Exp_\phi$ via Corollary \ref{C:Exp coeff}.
\end{proof}

\begin{proposition}
For $h\in N_\Gamma^*$, we have $J(h) = \delta_*^N(h)$.
\end{proposition}

\begin{proof}
For $w >0$, we express $h$ as in \eqref{E:h ell expression},
\[h^w  = d_{1,-w}^w \sigma^{-w}(h_1) + d_{2,-w}^w\sigma^{-w}(h_2) + \dots + d_{1,0}^w(h_1)+ \dots + d_{d,0}^w(h_d) + \mathcal O(\sigma).\]
Then we compute (recall that $\delta_0^N(\sigma^i(h_j))=0$ for all $i\geq 1$ and all $1\leq j\leq d$ and that $\delta_0(\mathcal O(\sigma))=0$)
\begin{align*}
J_k(h) &= \sum_{i=0}^\infty \left(\delta_0^M(\tau^{-i}(g_k))^\top\right )\twisti\delta_0^N(\sigma^{i}(h))\twisti\\
& = \sum_{i=0}^\infty \left(\delta_0^M(\tau^{-i}(g_k))^\top\right )\twisti  \delta_0^N\left(\sigma^{i}\left(\sum_{m=-w}^{k-w} \sum_{j=1}^d d^w_{j,m}\sigma^{m}(h_j) + \mathcal O(\sigma)\right )\right )\twisti\\
&=\sum_{i=0}^w \sum_{m=-w}^{-i} \sum_{j=1}^d \left(\delta_0^M(\tau^{-i}(g_k))^\top\right )\twisti d^w_{j,m}\delta_0^N(\sigma^{i+m}(h_j))\twisti     + \sum_{i=w+1}^\infty \left(\delta_0^M(\tau^{-i}(g_k))^\top\right )\twisti\delta_0^N(\sigma^{i}(h))\twisti\\
&=\sum_{m=-w}^0  \sum_{j=1}^d d^w_{j,m}\left[\sum_{i=0}^m\left(\delta_0^M(\tau^{-i}(g_k))^\top\right )\twisti \delta_0^N(\sigma^{i-m}(h_j))\twisti \right ]    + \sum_{i=w+1}^\infty \left(\delta_0^M(\tau^{-i}(g_k))^\top\right )\twisti\delta_0^N(\sigma^{i}(h))\twisti\\
&=\sum_{m=-w}^0  \sum_{j=1}^d d^w_{j,m}I_m(1,1,g_k,h_j)    + \sum_{i=w+1}^\infty \left(\delta_0^M(\tau^{-i}(g_k))^\top\right )\twisti\delta_0^N(\sigma^{i}(h))\twisti\\
&= d^w_{k,0} + \sum_{i=w+1}^\infty \left(\delta_0^M(\tau^{-i}(g_k))^\top\right )\twisti\delta_0^N(\sigma^{i}(h))\twisti.
\end{align*}
Note that in the last line in the above calculation, we have applied Proposition \ref{P:I vanishing}. We recognize 
\[\sum_{i=w+1}^\infty \left(\delta_0^M(\tau^{-i}(g_k))^\top\right )\twisti\delta_0^N(\sigma^{i}(h))\twisti\]
as the tail of $J_k(h)$, which is a convergent series. Thus taking the limit as $w\to \infty$ of both sides of the above equality causes the series tail to vanish and proves the theorem by Definition \ref{D:delta_* def}.
\end{proof}

\begin{theorem}\label{T:Main Log Alg Theorem}
Let $h\in N_\Gamma$ be such that 
\[h - \sigma(h) = \sigma(g),\]
for some $g\in N$. Then
\[\Exp_\phi(\delta_0^N(h)) = \delta_*^N(h).\]
In particular, if $\delta_*^N(h)$ is in $\overline K$, then $\delta_0^N(h)$ is log-algebraic.
\end{theorem}

\begin{proof}
Note that the functional equation for $h$ gives
\[h= \sigma^i(h) + \sum_{j=1}^i \sigma^j(g).\]
We then begin with Theorem \ref{T:Exponential from Pairing}. We recall that $\delta_1^N(\sigma(w)) = \delta_1^N(w)$ for all $w\in N$, that $\delta_1^N(h_k) = \be_k$ (the $k$th standard basis element), and that $\delta_0^N(\sigma^i(g)) = 0$ for $g\in N$. Using this, we find
\begin{align*}
\Exp_\phi(\delta_0^N(h)) &= \delta_1^N\left(\sum_{i=0}^\infty \sum_{k = 1}^{d} \delta_0^M(\tau^{-i}(g_k))^\top \delta_0^N(h) \sigma^i(h_k)  \right )\\
&= \sum_{i=0}^\infty \sum_{k = 1}^{d} (\delta_0^M(\tau^{-i}(g_k))^\top)\twisti \delta_0^N(h)\twisti \delta_1^N\left(\sigma^i(h_k)  \right )\\
&= \sum_{i=0}^\infty \sum_{k = 1}^{d} (\delta_0^M(\tau^{-i}(g_k))^\top)\twisti \delta_0^N\left(\sigma^i(h) + \sum_{j=1}^i \sigma^j(g)\right )\twisti \delta_1^N\left(h_k  \right )\\
&= \sum_{i=0}^\infty \sum_{k = 1}^{d} (\delta_0^M(\tau^{-i}(g_k))^\top)\twisti \delta_0^N(\sigma^i(h))\twisti \be_k\\
&= \sum_{k = 1}^{d}J_k(h)\be_k\\
&= \delta_*^N(h) ,
\end{align*}
where $\be_k$ is the $k$th standard basis vector.
\end{proof}

\subsection{Example: Carlitz Tensor Powers}\label{S:Catlitz Tensor Power Example}
In this brief subsection, we illustrate Theorem \ref{T:Main Log Alg Theorem} in the particular case of tensor powers of the Carlitz. In this simple case, we are able to achieve more specific formulas which leads to a pair of open questions at the end of this section. We continue with the notation and example of Example \ref{E:Tensor Power of Carlitz}. Let $h \in N_\Gamma$ such that
\[h-\sigma(h) = \sigma(g),\]
for some $g\in \C_\infty[t]$. This functional equation shows that $h$ has a meromorphic continuation to $\C_\infty$ with possible poles at $\theta^{q^i}$ for $i\geq 1$ of order at most $n$, and thus we are able to take residues of $hdt$ at the values $\theta^{q^i}$. Our Theorem \ref{T:Main Log Alg Theorem} then gives
\begin{align*}
\Exp_C\on(\delta_0^N(h)) &= \delta_1^N\left(\sum_{i=0}^\infty \sum_{k = 1}^{d} \delta_0^M(\tau^{-i}(g_k))^\top \delta_0^N(h) \sigma^i(h_k)  \right )\\
&= \delta_*^N\left(h\right )\\
\end{align*}
However, in this very concrete, basic case of tensor powers of Carlitz, we get more specialized formulas. Combining this all together, and applying \eqref{E:product rule} we get
\begin{align*}
\Exp_C\on(\delta_0^N(h)) &= \delta_1^N\left(\sum_{i=0}^\infty \sum_{k = 1}^{d} \delta_0^M(\tau^{-i}(g_k))^\top \delta_0^N(h) \sigma^i(h_k)  \right )\\
&= \delta_1^N\left(\sum_{i=0}^\infty \sum_{k = 1}^{d} \delta_0^M(\tau^{-i}(g_k))^\top \delta_0^N(\sigma^i(h)) \sigma^i(h_k)  \right )\\
&= \delta_1^N\left(\sum_{i=0}^\infty \sum_{k = 1}^{d} \delta_0^M(\tau^{-i}(g_k))^\top \delta_0^N((D_i(t)\twistk{1-i})^n h\twistk{-i}) \sigma^i(h_k)  \right )\\
&= \sum_{i\geq 0}\left (\partial_t^{n-1}\frac{(t-\theta)^{n-1}(D_i(t)\twist)^n h}{D_i(t)^n},\dots,\partial_t^{n-1}\frac{(D_i(t)\twist)^n h}{D_i(t)^n} \right )^\top\Bigg|_{t=\theta^{q^i}}\\
&= \sum_{i\geq 0}\left (\partial_t^{n-1}\left ((t-\theta)^{n-1}(t-\theta^{q^i})^n h\right ),\dots,\partial_t^{n-1}\left ((t-\theta^{q^i})^n h\right ) \right )^\top\Bigg|_{t=\theta^{q^i}}\\
&= \sum_{i\geq 0}\begin{pmatrix}
\res_{t=\theta^{q^i}} h dt\\
\res_{t=\theta^{q^i}} (t-\theta)h dt\\
\vdots\\
\res_{t=\theta^{q^i}} (t-\theta)^{n-1}h dt\\
\end{pmatrix}.
\end{align*}
Finally, using the fact that $h$ is regular away from $\theta^{q^i}$ for $i\geq 1$ and the fact that the sum of residues over all points of a curve must be zero, we conclude that the above sum equals
\[\Exp_C\on(\delta_0^N(h)) = 
\begin{pmatrix}
-\res_{\infty} h dt\\
-\res_{\infty} (t-\theta)h dt\\
\vdots\\
-\res_{\infty} (t-\theta)^{n-1}h dt\\
\end{pmatrix}.
\]
Thus in this case we see that $\delta_*^N(h)$ is related to the negative of the residue of $h$ at $\infty$. There are some hints towards a theory of this sort in Sinha's Main Diagram \cite[\S 4.2.3]{Sin97}. This invites the following pair of question related to residues.

\begin{question}
For a general Anderson $A$-module $\phi$, and $h \in N_\Gamma$ as above, can we relate $\Exp_\phi(\delta_0^N(h))$ to the residue of $h$ at $\infty$?
\end{question}

\begin{question}
In \cite[Thm. 2.5.13]{HJ20} the authors describe an isomorphism $\Xi$ from the algebraic dual of $M$ to $N$ constructed by taking a finite sum of residues at $\infty$. They pose an open question about a pairing derived from this map. How do the pairings $F$ or $G$ from our paper compare the map $\Xi$ from \cite{HJ20} and can we use them to construct this pairing explicitly?
\end{question}

\section{Acknowledgments}
The author give sincere thanks to the referee for a very careful reading of the paper and many constructive suggestions. He also thanks Andreas Maurischat for many helpful discussions about the context of the formulas contained in this paper.

%%%%%%%%%%%%%%%%%%%%%%%%%%%%%

%%%%%%%%%%%%%%%%%%%%%%%%%%%%%

\end{document}